\documentclass[12pt,reqno]{amsart}
\usepackage{epsfig}

\newcommand{\p}{\mathfrak{p}}
\newcommand{\q}{\mathfrak{q}}
\newcommand{\bigq}{\mathfrak{Q}}
\newcommand{\B}{\mathfrak{P}}
\newcommand{\m}{\mathfrak{m}}
\newcommand{\n}{\mathfrak{n}}
\newcommand{\aok}{\mathfrak{a}}

\newcommand{\OK}{{\mathcal{O}_{K}}}
\newcommand{\OL}{\mathcal{O}}
\newcommand{\F}{\mathcal{F}}
\newcommand{\G}{\mathcal{G}}
\newcommand{\E}{\mathcal{E}}
\newcommand{\M}{\mathcal{M}}
\newcommand{\SOP}{\mathfrak{E}} 
\newcommand{\pee}{\mathcal{P}}
\newcommand{\cond}{\mathfrak{f}}
\newcommand{\Q}{\mathbb{Q}}

\newcommand{\C}{\mathbb{C}}
\newcommand{\Z}{\mathbb{Z}}
\newcommand{\add}{{\rm add}}
\newcommand{\gal}{{\rm Gal}}
\newcommand{\ord}{{\rm ord}}
\newcommand{\End}{{\rm End}}

\newcommand{\diag}{{\rm diag}}
\renewcommand{\det}{{\rm det}}
\renewcommand{\char}{{\rm char}}

\theoremstyle{plain}
\newtheorem{thm}{Theorem}
\newtheorem{lem}{Lemma}[section]
\newtheorem{prop}[lem]{Proposition}
\newtheorem{cor}[lem]{Corollary}

\theoremstyle{definition}

\theoremstyle{remark}
\newtheorem*{rem}{Remark}
\newtheorem*{rems}{Remarks}
\newtheorem*{notation}{Notation}

\newsavebox{\proofbox}
\savebox{\proofbox}{\begin{picture}(7,7)%
  \put(0,0){\framebox(7,7){}}\end{picture}}

\begin{document}

\title{Coates-Wiles towers for CM abelian varieties}

\author[C. Rowe]{Christopher M. Rowe}
\address{Pacific Institute for the Mathematical Sciences, University of British Columbia, Room 205, 1933 West Mall, Vancouver, BC V6T 1Z2, Canada}
\email{rowec@math.ubc.ca}


\begin{abstract}
The aim of this paper is to compute congruence relations on units in fields generated by adjoining torsion points of a CM abelian variety to a number field.   For elliptic curves,  congruence relations of the form we compute were an important ingredient in the early proofs of the Coates-Wiles theorem.  In general, we compute congruence relations on exterior products of units in division fields, which more naturally fit into the framework of Rubin's generalization of Stark's conjecture.
\end{abstract} 

\maketitle

\section*{Introduction}

Let $E$ be an elliptic curve defined over $F$ with complex multiplication by the ring of integers of an imaginary quadratic extension $K$ of $\Q$, where $F$ is either $K$ or $\Q$.  In \cite{Coates-Wiles}, Coates and Wiles proved that if $E(K)$ has positive rank then the Hasse-Weil $L$-series $L(E/F,s)$ vanishes at $s=1$.  (In Rubin's important work on Tate-Shafarevich groups, he determined bounds on Selmer groups, which yield another proof of the Coates-Wiles theorem \cite{RubinTS}.)   

Coates and Wiles used formal groups and Iwasawa-theoretic techniques to relate elliptic units with special values of $L(E/F,s)$.  Using more classical techniques, Stark and Gupta were able to give a proof of the Coates-Wiles theorem for elliptic curves defined over $\Q$ \cite{StarkCW, Gupta}.  The proofs utilize different techniques, but both include many of the same ideas.

Specifically, let $E$ be an elliptic curve defined over $\Q$ with complex multiplication (CM) by the ring of integers of an imaginary quadratic extension $K$ of $\Q$ (so necessarily of class number one).  Let $\p = (\pi)$ be one of infinitely many suitably chosen primes of $K$, and $K_n$ the field of $\pi^n$-division values of $E$ (i.e., the finite extension of $K$ obtained by adjoining all the $\pi^n$-torsion of $E$ to $K$).  Then there exists a unique prime $\p_n$ of $K_n$ lying over $\p$.  Fix a point $Q \in E(\Q) \setminus \pi E(\Q)$ of infinite order and let $L_n$ be the field of $\p^n$-division values of $Q$, i.e., $L_n = K_n(\frac{1}{\p^n}Q)$.  
 
Gupta showed there exists a positive integer $e$ such that the
conductor of the field extension $L_e/K_e$ is $\p_e^2$, and then by class field theory all units of $K_e$ congruent to $1 \bmod \p_e$ are also congruent to $1 \bmod \p_e^2$.  Both Coates and Wiles and Gupta used a congruence relation of this form to show that $\pi | L(E/\Q,1)$, and hence $L(E/\Q,1)$ vanishes.

The goal of this paper is to compute conductors and congruence relations on units for CM abelian varieties paralleling the work of Gupta for CM elliptic curves.  In particular, let $A$ be an abelian variety of dimension $g$ defined over a Galois number field $K$ with complex multiplication by the ring of integers of $K$ (with the degree $2g$ over $\Q$).  Furthermore let $\p$ be a prime of $K$ lying over the odd rational prime $p$, and assume that both $p$ splits completely in $K$ and $A$ has good reduction at all primes lying over $p$.  Then we describe congruence relations on units of $K_n$, the extension of $K$ generated by the $\p^n$-torsion points of $A$.  For abelian surfaces, Grant described congruence relations on units similar to those of Gupta (see \cite{GrantCW}), but additional progress was stymied by the lack of understanding and construction of ``abelian units'' to parallel the theory of elliptic units.  Moreover, it seems to be a difficult problem to come up with a general theory of abelian units.

Since we require that $p$ splits completely, the field extensions $K_n/K$ can be shown to have degree $p^{n-1}(p-1)$ and to be totally ramified at half of the primes of $K$ lying above $p$ and unramified at the other half (which primes ramify depends upon $\p$ and the ``CM type'' of $A$).  We fix a point $Q \in A(K) \setminus \p A(K)$ of infinite order and construct field extensions $L_n = K_n(\frac{1}{\p^n}Q)$, the field of $\p^n$-division values of $Q$.  We are able to show that there exists a positive integer $e$ such that $L_n/K_n$ is unramified for $1 \leq n < e$, but $L_e/K_e$ is ramified.  Then using properties of formal groups, we are able to compute the conductor of $L_e/K_e$.

\begin{thm} \label{thm:cond}
There exists a set $\SOP$ of at least one and at most $g$ primes of $K_e$ lying over $p$ such that the conductor of $L_e/K_e$ is $\prod_{\B_e \in \SOP} \B_e^2$. 
\end{thm}

If $\SOP$ is composed of a lone prime, then we have computed a conductor exactly as in the work of Gupta and of Grant \cite{Gupta, GrantCW}.  In this case, we will have the same type of congruence relations on units of $K_e$.  If on the other hand $\SOP$ is composed of more than one prime, say $\# (\SOP) = s > 1$, we have a congruence relation on exterior products of units of $K_e$.  

\begin{thm} \label{thm:cong}
Let $u_1, \ldots, u_s$ be units of $K_e$, which are
congruent to $1 \bmod \prod_{\B_e \in \SOP} \B_e$.   Then $u_1 \wedge \cdots
\wedge u_s$ is trivial $\bmod$  $\prod_{\B_e \in \SOP} \B_e^2$.
\end{thm}

This fits more naturally into the framework of Rubin's generalization of Stark's conjecture \cite{RubinS, StarkIV}.  Stark's conjecture is a generalization of the class number formula, which relates the arithmetic of number fields to special values of Artin $L$-series.  For example, let $F/E$ be a finite abelian extension of number fields and $\chi$ a character on  $G = \gal(F/E)$.  Let $S$ be a finite set of primes of $E$ including the archimedian primes, the primes which ramify in $F$, and a set of $r$ primes which split completely in $F$.  Stark's conjecture relates the lead term of the Taylor expansion of the Artin $L$-series $L(s,\chi)$ at $s=0$ to the determinant of an $r \times r$ matrix whose entries are linear combinations of logarithms of absolute values of $S$-units in $F$ (i.e., units locally at all primes not in $S$).  For $r=1$, Stark gave the following refined conjecture.  Let $\B_1$ be any prime of $F$ sitting over the designated prime of $S$.  Then there exists an $S$-unit $\epsilon_1$ such that for all characters $\chi$ on $G$
\begin{equation*}
L'(0,\chi) = -\frac{1}{w_F}\sum_{\gamma \in G} \chi(\gamma) \log |\gamma(\epsilon_1)|_{\B_1}. 
\end{equation*}

When $F = \Q$ or an imaginary quadratic extension of $\Q$, Stark was able to use properties of cyclotomic and elliptic units, respectively, to prove his refined conjecture \cite{StarkIV}.  For $r \geq 1$, Rubin gave a generalized refined Stark's conjecture, which conjectured a relation between exterior products of $S$-units and the lead term in the Taylor expansion of  $L(s,\chi)$ at $s=0$ \cite{RubinS}.  Rubin's conjecture relates an entry in a lattice of $\Q \otimes \bigwedge^r U_{S,T}$ to the $r$th derivative at $s=0$ of $L(s, \chi)$, where $U_{S,T}$ is a specific subgroup of the group of $S$-units, depending upon an auxiliary, finite set of primes $T$.    

Rubin's conjecture applied to $K_e/K$ should produce exterior products of ``abelian $S$-units'' in $K_e$.  Outside of Rubin's original paper, the only direct evidence for Rubin's conjecture is given in \cite{GrantSR}, which looks at exterior products of units arising from $5$-torsion on the Jacobian of $y^2=x^5 + 1/4$.  However, for $g=1$, Stark further refined his conjecture so that the $S$-unit in $K_n$ relating to the $L$-series is actually a unit, and it was congruences on units, not the corresponding $S$-units, that were employed in the proof of the Coates-Wiles theorem.  These results led us to consider whether the existence of a point of infinite order in $A(K)$ should force congruence conditions on \emph{exterior products} of units in $K_e$; the result of which is Theorem \ref{thm:cong}.

The first section of this paper fixes notation and assumptions about our abelian variety $A$ and number field $K$.  In the second section, we collect the information we need about formal groups attached to abelian varieties.  Then we use properties of formal groups attached to abelian varieties in sections \ref{sec:divfld} and \ref{sec:tower} to describe properties of the field extensions $K_n/K$ and $L_n/K_n$, respectively.  We prove Theorem \ref{thm:cond}  in section \ref{sec:cond} and Theorem \ref{thm:cong} in section \ref{sec:units}.

Most of the results of this paper were contained in the author's Ph.D. thesis, and the author would be remiss if he did not thank his advisor, David Grant, for his invaluable assistance.  Also, the author would like to thank both Wolfgang Schmidt and the Pacific Institute for the Mathematical Sciences for their support during the writing of this paper.

\begin{notation}  
Let $F/E$ be number fields, and $\q$ a prime of $E$.  We let $\cond(F/E)$ and $D(F/E)$ denote the conductor and discriminant of $F/E$ respectively.  We let $\overline{E}$ denote an algebraic closure of $E$, $\OL_E$ the ring of integers of $E$, $E_{\q}$ the completion of $E$ at $\q$, and $\OL_{\q}$ the ring of integers of $E_{\q}$.  
\end{notation}

\section{The Setup}\label{sec:setup}

We will use the terminology and results of the theory of complex multiplication as developed by Lang in \cite{LangCM} for the results in this section. 

Throughout this paper, $K$ will denote a number field of degree $2g$ over $\Q$, which is both Galois and a CM field.  Recall that $K$ is a CM field if it is a totally imaginary quadratic extension of a totally real number field.  Furthermore we assume that $A$ is an abelian variety of dimension $g$ with complex multiplication by the ring of integers of $K$.  This means that we have an embedding
\begin{equation*}
i:K \hookrightarrow \End(A) \otimes \Q = \End_{\Q}(A)
\end{equation*}  
\noindent
such that $i(\OK) = i(K) \cap \End(A)$.  

The action of $i$ on the tangent space of $A$ determines a CM type $\Phi = \{ \phi_1, \ldots, \phi_g \}$, where $\phi_i, \phi_j$ are non-conjugate embeddings of $K$ into $\C$ for $1 \leq i,j \leq g$.

A CM abelian variety will be a pair $(A,i)$.  We say that $(A,i)$ is defined over a number field $K$ if both $A$ and every element of $\End(A)$ are defined over $K$.  Let $\overline{K}$ be a fixed algebraic closure of $K$.  For $\alpha \in \OK$, we put $[\alpha] = i(\alpha): A(\overline{K}) \rightarrow A(\overline{K})$.  Therefore we have that $\sigma ([\alpha]) = [\sigma (\alpha)]$ for any $\sigma \in \gal(\overline{K}/K)$.   We let $A[\alpha] = \mbox{ker}[\alpha]$ denote the $\alpha$-torsion of $A$.



Throughout this paper $p$ will represent an odd rational prime such that (i) $A$ has good reduction at all primes lying over $p$ and (ii) $p$ splits completely in $K$.  Furthermore we fix a prime $\p$ of $K$ such that $\p|p$ and an element $\pi \in \OK$ such that $\mbox{ord}_p \pi = \ord_{\p} \pi = 1$ (such a $\pi$ exists by the Chinese remainder theorem since $p$ splits completely).  For ease of notation, we put $\p_i = \phi_i(\p)$ for $1 \leq i \leq g$ and let $S = \{\p_1, \ldots, \p_g \}$.

\begin{rem}
The Jacobians of rational images of Fermat curves are a class of abelian varieties satisfying the assumptions of this section  (see \cite{LangCM} or \cite{Shimura} for more details).
\end{rem}

\section{Formal groups}\label{sec:fmlgrp}

We need a variety of results on formal groups.  We refer to Hazewinkel's excellent book \cite{Hazewinkel} as a general reference on formal groups and to Hindry and Silverman \cite{Hindry-Silverman} for the construction of a formal group attached to an abelian variety and results on the kernel of reduction of $A$ mod $\p$.

\subsection{Basic Properties}

Let $R$ be a commutative ring with identity and let $X = (X_1, \ldots, X_n)^t$, $Y = (Y_1, \ldots, Y_n)^t$, and $Z = (Z_1, \ldots, Z_n)^t$ be column vectors of variables.  We call an $n$-tuple of power series over $R$, $F = (F_1, \ldots, F_n)^t$ in $2n$-variables an $n$-dimensional \emph{commutative formal group (law)}  over $R$ if
\begin{align*}
  F(X,Y) &= X + Y + (d^{\circ} \geq 2), \\
 F(X,F(Y,Z)) &= F(F(X,Y),Z) \mbox{, and} \\
F(X,Y) &= F(Y,X),
\end{align*}
\noindent
where $(d^{\circ} \geq m)$ denotes a $n$-tuple of power series, all of whose terms are of total degree at least $m$.

If $F$ and $G$ are two $n$-dimensional formal groups defined over $R$.  A \emph{homomorphism} $\varphi: F \rightarrow G$ over $R$ is a $n$-tuple of power series over $R$ without constant terms such that $\varphi (F(X,Y)) = G(\varphi(X), \varphi(Y))$.  A homomorphism is an \emph{isomorphism} if it has a two-sided inverse.  Let $\varphi = (\varphi_1, \ldots, \varphi_n)^t$ be a homomorphism from $F$ to $G$.  Suppose that the linear term of $\varphi_i(X_1, \ldots, X_n)^t$ is $\sum_{j=1}^n a_{ij}X_j$.  We call the matrix $(a_{ij})$ the \emph{jacobian} of $\varphi$, which we denote by $j(\varphi)$.  The following is elementary.  

\begin{lem} \label{fmlgrpjac}
A homomorphism $\varphi$ between two $n$-dimensional formal groups defined over $R$ is an isomorphism if and only if $\det(j(\varphi))$ is a unit in $R$.   
\end{lem}

Furthermore, we call two formal groups \emph{strictly isomorphic} over $R$ if there is an isomorphism $\varphi$ over $R$ such that $j(\varphi)$ is the identity matrix.  In \cite{Hindry-Silverman}, Hindry and Silverman demonstrate how to construct a $g$-dimensional commutative formal group $\F$, defined over $K$, associated to our $g$-dimensional abelian variety $A$ defined over $K$.  Moreover, $\F$ can be defined over $\Z_p$ and, since $p$ splits completely, we will show that $\F$ is strictly isomorphic to a product of $g$ one-dimensional formal groups defined over $\Z_p$ of Lubin-Tate type.  In order to do this, we need to formalize what we mean by a ``product of formal groups''.

Let $F=(F_1, \ldots, F_n)^t, G=(G_1, \ldots, G_m)^t$ be $n$ and $m$-dimensional formal groups, respectively, defined over $R$.  Let $X=(X_1, \ldots, X_n)^t$ and $Y=(Y_1, \ldots, Y_n)^t$ be $n$-tuples of variables such that $F_i(X,Y) \in R[[X,Y]]$ for $1 \leq i \leq n$.  Also let $W = (W_1, \ldots, W_m)^t$ and $Z=(Z_1, \ldots, Z_m)^t$ be $m$-tuples of variables such that $G_i(W,Z) \in R[[W,Z]]$ for $1 \leq i \leq m$.  Then we can define an $(m +n)$-dimensional formal group $H$ over $R$ with $H_i \in R[[X,W,Y,Z]]$ by putting 
\begin{align*}
H_i & = F_i \mbox{ for $1 \leq i \leq n$ and } \\
H_{n+i} & = G_i \mbox{ for $1 \leq i \leq m$.}
\end{align*}

Specifically, if $\G_1, \ldots, \G_g$ are one-dimensional commutative formal groups over $R$, then iterating the above construction gives us a $g$-dimensional commutative formal group over $R$.  Let $\G$ be the formal group constructed in this manner, then we will call such a $\G$ a \emph{product of one-dimensional commutative formal groups}, and write $\G = \oplus_{j=1}^g \G_j$.  Let $\varphi_j$ be an endomorphism of $\G_j$ for $1 \leq j \leq g$.  Then it is easy to see that $\varphi = (\varphi_1, \ldots, \varphi_g)^t$ is an endomorphism of $\G$, which we will denote by $\varphi = \oplus_{j=1}^g \varphi_j$.

\subsection{Properties of the Formal Group $\F$ Associated with $A$}

Following Hindry and Silverman \cite{Hindry-Silverman}, we associate to $A$ a $g$-dimensional commutative formal group $\F$ over $K$.  Let $O$ denote the origin of $A$, $\widehat{\OL}_A$ the completed local ring at the origin of $A$, and $\n$ the maximal ideal of the local ring at the origin of $A$.  Since $A$ is nonsingular, there is an isomorphism $\displaystyle \widehat{\OL}_A \cong K[[s_1, \ldots, s_g]]$, where $s_1, \ldots, s_g \in \n$ are fixed local parameters on $A$ at the origin.   

Next we consider the product $A \times A$.  Let $p_i:A \times A \rightarrow A$ be the projection onto the $i$th factor for $i = 1,2$ and for local parameters at the point $(O,O) \in A \times A$ we choose the functions $x_1, \ldots, x_g, y_1, \ldots, y_g$, where $x_i:= s_i \circ p_1$ and $y_i:= s_i \circ p_2$.  Let $\widehat{\OL}_{A \times A}$ be the completed local ring at the origin of $A \times A$.  Then this choice determines an isomorphism $\widehat{\OL}_{A \times A} \cong \widehat{\OL}_A \times \widehat{\OL}_A \cong  k[[x_1, \ldots, x_g, y_1, \ldots, y_g]]$.  

The morphism giving the group law on $A$, $\add:A \times A \rightarrow A$, induces a map of local rings $\add^*: \widehat{\OL}_A \rightarrow \widehat{\OL}_{A \times A}$  and by the above isomorphisms, a map $\add^*: k[[s_1, \ldots, s_g]] \rightarrow k[[x_1, \ldots, x_g, y_1, \ldots, y_g]]$ of formal power series rings.  We let $\F_i= \add^*(s_i)$, and $\F = (\F_1, \dots, \F_g)^t$.  However the formal group $\F$ depends upon the choice of parameters $s_1, \ldots, s_g \in \n$.  Since $A$ is defined over $K$ and the determinant of a change of basis matrix of $\n/\n^2$ is a unit in $K$, Lemma \ref{fmlgrpjac} shows that all formal groups associated to $A$ in this way are isomorphic over $K$.  Since $A$ has complex multiplication by $\OK$ and $K$ is Galois over $\Q$, we are able to show that the formal group $\F$ possesses a ``CM action''.

\begin{lem} \label{Kcmact}
There exist local parameters $s_1, \ldots, s_g$ at the origin of $A$ defined over $K$ such that for any $\alpha \in \OK$ 
\begin{equation} \label{cmact}
[\alpha]^* s_i := s_i \circ [\alpha] = \phi_i(\alpha)s_i + (d^{\circ} \geq 2) \mbox{, in $\widehat{\OL}_A$}.
\end{equation}
\end{lem}

\begin{proof}
Let $\B$ be a prime of $K$ lying over $p$, then $A$ has good reduction at $\B$ by assumption.  Let $\tilde{A}$ denote the reduced abelian variety $\bmod$ $\B$, which is defined over $\tilde{K} = \OK/ \B \OK$, and $H^0(A, \Omega)$ the space of holomorphic differentials of $A$.  
There exists a basis $\omega_1, \ldots, \omega_g$ of differentials defined over $K$ such that $[\alpha]^* \omega_i := \omega_i \circ [\alpha] = \phi_i(\alpha) \omega_i$ for all $1 \leq i \leq g$ and $\alpha \in \OK$.  Indeed, this follows from \cite[pg. 99]{Shimura}, since $K/\Q$ is Galois, the field of definition of $(A,i)$ is $K$, and $(A,i)$ has CM type $(K, \Phi)$.  We will call such a basis a \emph{splitting basis} for $H^0(A, \Omega)$. 

Since $\omega \in H^0(A, \Omega)$ is translation invariant, we can characterize $\omega$ by its representation as a differential at the origin \cite[pg. 13]{Shimura}.  Recall that $\n$ is the maximal ideal of the local ring at the origin of $A$, and let $\psi$ be the isomorphism $\n/\n^2 \cong H^0(A, \Omega)$ (see \cite{Milne}), where $\tau \in \n$ maps to the differential at the origin represented by $d\tau$.  Now we pick a $S_i \in \n$ such that $\psi(S_i) = \omega_i$, and since $\psi$ commutes with endomorphisms of $A$ \cite[pg. 75]{Shimura}, we have that $[\alpha]^*S_i = \phi_i(\alpha)S_i + (d^{\circ} \geq 2)$.  

We still need to show that we can choose parameters defined over $K$.  Since $\gal(\overline{K}/K)$ commutes with $\psi$ (see \cite{Shimura}), we see that $\sigma(S_i) = S_i \bmod \n^2$ for all $\sigma \in \gal(\overline{K}/K)$.  If we let $E = K(S_i)$, then $E/K$ is a finite Galois extension.  Now we set $s_i = (1/|\gal(E/K)|)\sum_{\tau \in \scriptstyle{\gal}(E/K)}\tau(S_i)$, which gives us the desired parameters, now defined over $K$.
\end{proof}

By Lemma \ref{Kcmact}, we can define $([\alpha]_{\F})_i(s_1, \ldots, s_g) := [\alpha]^*s_i \in K[[s_1, \ldots, s_g]]$ for any $\alpha \in \OK$.  Then the map $\alpha \mapsto [\alpha]_{\F} = (([\alpha]_{\F})_1, \ldots, ([\alpha]_{\F})_g)^t$ gives an embedding of $\OK$ into the endomorphism ring of $\F$.  Next we show that $\F$ and $[\alpha]_{\F}$ are actually defined over $\OK$. 

\begin{lem} \label{Kpcmact}
There exist local paramaters $s_1, \ldots, s_g$ at the origin of $A$ such that $\F$ is defined over $\OK$ and $[\alpha]_{\F}$ is a $g$-tuple of power series in $\OK[[s_1, \ldots, s_g]]$ for every $\alpha \in \OK$.  Moreover, if $s = (s_1, \ldots, s_g)^t$ then    
\begin{equation*}
([\alpha]_{\F} s)_i = \phi_i(\alpha) s_i + (d^{\circ} \geq 2).
\end{equation*} 
\end{lem}

\begin{proof}
Let $\B$ be a prime of $K$ lying over $p$.  
In the proof of Lemma \ref{Kcmact}, we saw that there exists a splitting basis $\omega_1, \ldots, \omega_g$ for $H^0(A, \Omega)$, where the $\omega_i$ are defined over $K$.    

For $\omega \in H^0(A, \Omega)$, a suitable multiple of $\omega$ will reduce to a non-zero differential on $\tilde{A}$ \cite[pg. 80]{Shimura}.  When $\omega$ reduces to a non-zero differential on $\tilde{A}$, we let $\widetilde{\omega}$ denote the reduced differential.  
After multiplication by a suitable multiple $a_i \in K$, $\widetilde{a_i\omega_i}$ is non-zero and defined.  If we rename $a_i \omega_i$ as $\omega_i$, then the (new) $\omega_i$ form a basis for $H^0(A, \Omega)$ such that $[\alpha]^* \omega_i = \phi_i(\alpha) \omega_i$, and hence the CM action is preserved.  

Let $\psi$ be the isomorphism between $\n/\n^2$ and $H^0(A, \Omega)$ given in the proof of Lemma \ref{Kcmact} and $s_1, \ldots, s_g$ be local parameters at the origin of $A$ over $K$ such that $\psi(s_i) = \omega_i$.  Since $p$ is unramified in $K$ by assumption, the $\widetilde{\omega_i}$ form a basis for $H^0(\tilde{A}, \Omega)$ \cite[pg. 99]{Shimura}.  
Since $\psi$ commutes with the reduction map $\bmod$ $\B$, letting $\widetilde{s_i}$ be the reduction of $s_i$, the $\widetilde{s_i}$ form a set of parameters at the origin of $\tilde{A}$.  
Since $A$ has good reduction at $\B$, the addition map on $A$ reduces to the addition map on $\tilde{A}$.  Therefore we can construct a formal group on the reduced abelian variety coming from the parameters $\widetilde{s_1}, \ldots, \widetilde{s_g}$, which will yield power series $F_i \in (\OK/\B \OK)[[X_1, \ldots, X_g, Y_1, \ldots, Y_g]]$.  
By construction, $\widetilde{add^{*}s_i}= \widetilde{add^{*}}\widetilde{s_i}$, 
so the formal group $\F$ defined by $s_1, \ldots, s_g$ over $K$ reduces $\bmod$ $\B$ to  $\widetilde{\F_i} = F_i$.  Therefore $\F_i \in \OK[[X_1, \ldots, X_g, Y_1, \ldots, Y_g]]$.   

Since endomorphisms of $A$ reduce to endomorphisms of $\tilde{A}$, endomorphisms of $\F$ reduce to endomorphisms of $\tilde{\F}$.  Indeed, $\widetilde{[\alpha]^{*}s_i} = \widetilde{[\alpha]}^{*}\widetilde{s_i}$ (see \cite[pg. 75]{Shimura}), and so we must have $[\alpha]^*s_i \in \OK[[s_1, \ldots, s_g]]$. 
\end{proof}

\begin{rems} ~
\begin{enumerate}
\item For a CM elliptic curve, the CM type is usually taken to consist of the identity, so the CM action is trivial.
\item We will use this CM action to determine which primes ramify.  (We became aware of this use of the CM action in \cite{Grant5}.)
\end{enumerate}
\end{rems}

Let $\Phi_K^{'} = \{ \phi_1^{-1}, \ldots, \phi_g^{-1} \}$, then $(K,\Phi_K^{'})$ is a CM type \cite[pg. 62]{LangCM}.  Let $N_{\Phi_K^{'}}(x) = \prod_{\phi \in \Phi_K^{'}} \phi(x)$ for $x \in K$ and extend $N_{\Phi_K^{'}}$ this to ideals of $K$ in the usual way.  Let $\B$ be a prime of $K$ lying over $p$, then $A$ has good reduction at $\B$ by assumption.  Since $A$ has principal complex multiplication by $K$, there exists an element $\alpha_{\B} \in \OK$ such that $[\alpha_{\B}]$ reduces to the Frobenius endomorphism $\bmod$ $\B$ (i.e., the endomorphism $x \mapsto x^{N\B}$ on $\tilde{A}$) \cite[pg. 61]{LangCM}.    

\begin{lem} \label{frobdecomp}
Let $\B$ be a prime of $K$ lying over $p$, and $\alpha_{\B}$ the element of $\OK$ which reduces to the Frobenius $\bmod$ $\B$.  Then $\alpha_{\B}$ has ideal decomposition in $K$ given by
\begin{equation}\label{frobenius}
(\alpha_{\B}) = N_{\Phi_K^{'}}(\B) = \prod_{j = 1}^g \phi_j^{-1}(\B). 
\end{equation}
\noindent
Moreover, $\ord_{\p} \alpha_{\B} = 1$ if and only if $\B$ is also in $S$.
\end{lem}

\begin{proof}
The left hand equality in $\eqref{frobenius}$ is just \cite[pg. 88]{LangCM} and the right hand equality comes from the definition of $N_{\Phi_K^{'}}$.  Since $p$ splits completely and $K/\Q$ is Galois, the $\phi_j^{-1}(\B)$ are distinct primes lying over $p$.  On the one hand, if $\B = \p_i$ for some $1 \leq i \leq g$, then $\phi_i^{-1}(\phi_i(\p)) = \p$ is a term in the product.  On the other hand, if $\B|p$, but $\B \notin S$, then $\p \neq \phi_j^{-1}(\B)$ for any $1 \leq j \leq g$. 
\end{proof}

\subsection{The Kernel of Reduction}

Fix a prime $\B$ of $K$ lying over $p$.  By base extension, we consider $\F$ defined over $\OL_{\B}$.    Let $L$ be a finite extension of $K_{\B}$, $\OL_{L}$ its ring of integers, and $\m$ its maximal ideal.  By assumption, $A$ considered over $\OL_L$ has good reduction at $\m$, so we let $\tilde{A}$ denotes the reduced abelian variety $\bmod$ $\m$, which is defined over $\tilde{L} = \OL_{L}/ \m \OL_{L}$.  We define the \emph{kernel of reduction} of $A \bmod \m$ by 
\begin{equation} \label{krdef}
A^{\circ}(L) := \mbox{ker} \left\{ A(L) \stackrel{red}{\longrightarrow}\tilde{A}(\tilde{L}) \right\}.
\end{equation}

Now consider $X,Y \in \m^g = \m \times \cdots \times \m$.  Then $\F_i(X,Y)$ will converge in $\OL_L$.  Indeed, $\OL_L$ is a complete local ring, and $\F$ is defined over $\OL_{\B} \subseteq \OL_L$.  Thus the formal group $\F$ defines a group structure on $\m^g$.  Let $\F(\m)$ be the set of $g$-tuples $\m^g$ with the group law
\begin{equation} \label{group}
\m^g \times \m^g \stackrel{+_{\F}}{\longrightarrow} \m^g \mbox{, given by} \ \  X +_{\F} Y := \F(X,Y).
\end{equation}
\noindent
This gives us the following isomorphism
\begin{equation} \label{kr}
A^{\circ}(L) \cong \F(\m),
\end{equation}
\noindent
given by $A^{\circ}(L) \ni R \mapsto (s_1(R), \ldots, s_g(R)) \in \m^g$ \cite[Thm. C.2.6]{Hindry-Silverman}.

Now let $\overline{K_{\B}}$ be a fixed algebraic closure of $K_{\B}$ with valuation ring $\OL$ and maximal ideal $\M$.  Although $\overline{K_{\B}}$ is not complete, each of its elements lives in a finite extension of $K_{\B}$.  So we extend $+_{\F}$ to $\M^g$ and identify $\F(\M)$ with $A^{\circ}(\overline{K_{\B}})$ (defined analogously to $\eqref{krdef}$).  By Lemma \ref{Kpcmact}, $[\alpha]_{\F}$ is an endomorphism of $\F$ for every $\alpha \in \OK$.  For $\alpha \in \OK$, we define the \emph{$\alpha$-torsion of $\F$} to be  $\F[\alpha] = \mbox{ker} \{[\alpha]_{\F}: \F(\M) \rightarrow \F(\M) \}$, and for any ideal $\aok \in \OK$, we define $\F[\aok] = \cap_{\alpha \in \aok} \F[\alpha]$ to be the \emph{$\aok$-torsion} of $\F$.  

\begin{lem} \label{primetop} ~
\begin{itemize}
\item[(i)] $\F(\M)$ has no torsion relatively prime to $p$.
\item[(ii)]  $A[\p^n]$ is in the kernel of reduction $\bmod$ $\M$ if and only if $\B \in S$.
\end{itemize}
\end{lem}

\begin{proof}
(i)  This is Proposition C.2.5 of \cite{Hindry-Silverman}.

(ii)  On the one hand, if $\B \notin S$ then $[\pi^n]_{\F}$ is an endomorphism of $\F$.  So by Lemma \ref{fmlgrpjac} the determinant of the jacobian of $[\pi^n]_{\F}$ is $\det(j([\pi]^n_{\F})) = \prod_{j = 1}^g \phi_j(\pi^n)$.  As in the proof of Lemma \ref{frobdecomp}, we see that $\ord_{\B} \det(j([\pi^n]_{\F})) = 0$ and therefore $j([\pi^n]_{\F})$ is invertible.  Hence $[\pi^n]_{\F}$ is an automorphism of $\F(\M)$, so $\F(\M)$ can have no $\pi^n$-torsion.  By definition, $\F[\p^n] \subseteq \F[\pi^n]$, and the result follows.

On the other hand, if $\B \in S$ and $\alpha_{\B}$ is the element of $\OK$ which reduces to the Frobenius $\bmod$ $\B$.  Now $[\alpha_{\B}^n]$ is a purely inseparable morphism, so $\tilde{A}[\alpha_{\B}^n] = \tilde{O}$.  Therefore $A[\alpha_{\B}^n]$ is in the kernel of reduction $\bmod$ $\B$.  Basic facts about torsion groups and Lemma \ref{frobdecomp} give that $A[\alpha_{\B}^n] = A[\p^n] \oplus A[\aok]$ for some integral ideal $\aok$ of $K$.  Hence $A[\p^n]$ is in the kernel of reduction $\bmod$ $\B$.
\end{proof}

\begin{cor} \label{AFtors}
Let $\B \in S$ and consider $\F$ defined over $\OL_{\B}$.  Then we can identify $A[\p^n]$ with $\F[\pi^n]$. 
\end{cor}

\begin{proof}
Lemma \ref{primetop} allows us to identify $A[\p^n]$ with $\F[\p^n]$, and we know that $\F[\p^n] \subseteq \F[\pi^n]$.  Also, the ideal $\p^n$ is generated by $\pi^n$ and $\p$ for $m > n$.  We take $m$ to be a multiple of the class number of $K$, so that $\p = (\gamma)$ for some $\gamma \in \OK$.  Then $\pi = \gamma \delta$ with $\delta$ relatively prime to $p$ by our choice of $\pi$, so $[\delta]$ is an automorphism of $\F$ over $\OL_{\B}$.  Therefore $\F[\p] = \F[\pi]$, and hence $\F[\p^n] = \F[\p] \cap \F[\pi^n] = \F[\pi] \cap \F[\pi^n] = \F[\pi^n]$.  
\end{proof}

\subsection{A Product of Formal Groups}

Now we will construct a formal group $\G$ strictly isomorphic to $\F$ over $\OL_{\B}$ for (fixed) $\B \in S$, where $\G$ is the product of one-dimensional Lubin-Tate formal groups.  This added structure will be very useful in what follows.  In order to show that $\F$ is isomorphic to a product of one-dimensional formal groups, we need to recall the properties of higher-dimensional Lubin-Tate formal groups (see \cite{Hazewinkel}).

\begin{rem}
In \cite{deShalit}, de Shalit gave a short proof of this decomposition using the theory of $p$-divisible groups.  However, we need to be careful to show that this isomorphism preserves the CM action as given by Lemma \ref{Kpcmact}, and we could not find this anywhere in the literature. 
\end{rem}

Let $\pi_{\B}$ be a prime element of $\OL_{\B}$, i.e., $(\pi_{\B}) = \B \OL_{\B}$.  Let $M$ be a $g \times g$ matrix such that $\pi_{\B}^{-1} M$ is an invertible matrix with entries in $\OL_{\B}$.  We let $\E_M$ denote the set of all $g$-tuples of power series $d(X)$ in $X= (X_1, \ldots, X_g)^t$ such that 
\begin{equation} \label{lubin-tate-cong}
d(X) = MX + (d^{\circ} \geq 2), \ \   \ \ d(X) \equiv X^p \bmod (\pi_{\B}).
\end{equation}
\noindent
Then the following is \cite[Thm. 13.3.3]{Hazewinkel} adapted to our assumptions.

\begin{lem} \label{Haz13.3.3}
For each $d(X) \in \E_M$, there is  precisely one $g$-dimensional formal group $F_d(X,Y)$ over $\OL_{\B}$ such that $F_d(d(X),d(Y)) = d(F_d(X,Y))$, so $d \in \End(F_d)$.  If $d(X), \overline{d}(X) \in \E_M$, then $F_d(X,Y)$ and $F_{\overline{d}}(X,Y)$ are strictly isomorphic over $\OL_{\B}$.
\end{lem}

For $g =1$, a formal group satisfying Lemma \ref{Haz13.3.3} is called a \emph{Lubin-Tate  formal group}, and hence is commutative \cite{Lubin-Tate}.  For $g>1$, we call a formal group $F$ with an endomorphism $d$ as in $\eqref{lubin-tate-cong}$ a \emph{higher dimensional Lubin-Tate formal group}.

For the rest of this section, fix an $i$ with $1 \leq i \leq g$.  Then we let $\B = \phi_i(\p)$,  $\alpha_{\B} \in \OK$ be such that $[\alpha_{\B}]$ reduces to the Frobenius $\bmod$ $\B$  and $s = (s_1, \ldots, s_g)^t$ represent local parameters over $K$ at the origin of $A$ which are parameters for a $g$-dimensional commutative formal group $\F$ over $\OL_{\B}$ with CM action as in Lemma \ref{Kpcmact}.  It will be useful to recall that $p$ splits completely and $\pi \in \OK$ such that $\ord_{\p} \pi = \ord_p \pi = 1$.

\begin{lem} \label{hdlt}
$\F$ is a higher dimensional Lubin-Tate formal group defined over $\OL_{\B} \cong \Z_p$.
\end{lem}

\begin{proof}
By Lemma \ref{Kpcmact}, we have that $[\alpha_{\B}]^* s_j = \phi_j(\alpha_{\B}) s_j + (d^{\circ} \geq 2)$.  But by Lemma \ref{frobdecomp} $(\alpha_{\B}) = \prod_{k=1}^g \phi_k^{-1}(\B)$, so $\B$ exactly divides $\phi_j(\alpha_{\B})$ for all $1 \leq j \leq g$.  
Let $M$ be the diagonal matrix whose $i$th diagonal entry is given by the coefficient of the linear term of $[\alpha_{\B}]^*s_j$ for each $1 \leq j \leq g$, i.e., $M = \diag\langle \phi_1(\alpha_{\B}), \ldots, \phi_g(\alpha_{\B}) \rangle$.  Let $\pi_{\B}$ be chosen by the Chinese Remainder Theorem such that $\B$ exactly divides $(\pi_{\B})$ and no other conjugate of $\B$ divides $(\pi_{\B})$.  Then we can write $\phi_j(\alpha_{\B}) = \pi_{\B} u_j$ with $u_j \in \OK$ and $u_j$ relatively prime to $\B$, and hence $M = \pi_{\B} \diag\langle u_1, \ldots, u_g \rangle$.  The $u_j$ are units in $\OL_{\B}$, so $\pi_{\B}^{-1} M$ is an invertible matrix.  
Then we have, by the definition of the Frobenius and since $\B$ is a first degree prime, 
\begin{equation*}
[\alpha_{\B}]_{\F}(s) = Ms + (d^{\circ} \geq 2), \ \  \ \  [\alpha_{\B}]_{\F}(s) \equiv s^p \bmod (\pi_{\B}).
\end{equation*}      
\noindent
Hence $[\alpha_{\B}]_{\F} \in \E_M$, and $\F$ is a higher dimensional Lubin-Tate formal group.
\end{proof}

Now we are in a position to decompose $\F$ into a product of one-dimensional formal groups over $\B$.  

\begin{prop} \label{fmlgrpdecomp}
Let everything be as in Lemma \ref{hdlt}.

\begin{itemize}
\item[(i)]  $\F$ is strictly isomorphic to a product of $g$ one-dimensional commutative formal groups of Lubin-Tate type defined over $\OL_{\B}$, say $\G = \oplus_{j=1}^g \G_j$, where $\G$ is given by paramaters $t= (t_1, \ldots, t_g)^t$.

\item[(ii)]  For each $1 \leq j \leq g$, there is an embedding of $\OK$ into $\End(\G_j)$ given by $\alpha \mapsto [\alpha]_{\G_j}$, where $[\alpha]_{\G_j} = \alpha t_j + (d^{\circ} \geq 2)(t_j)$.  Then we have an embedding $\OK$ into $\End(\G)$ given by $\alpha \mapsto \oplus_{j=1}^g [\phi_j(\alpha)]_{\G_j}$.  Moreover, we can identify $\G[\alpha]$ with $\F[\alpha]$ for all $\alpha \in \OK$.  

\item[(iii)]  Let $\pi_i = \phi_i(\pi)$, then $\ord_{\B} \pi_i = 1$ and $A[\p^n]$ can be identified with $\G_i[\pi_i^n]$.  
\end{itemize} 
\end{prop}

\begin{proof}
(i)  Let $M = \diag \langle \phi_1(\alpha_{\B}), \ldots, \phi_g(\alpha_{\B}) \rangle$ and $\pi_{\B}$ be as in the proof of Lemma \ref{hdlt}, then $M = \pi_{\B} \diag \langle u_1, \ldots, u_g \rangle$, where $[\alpha_{\B}]_{\F}s = Ms + (d^{\circ} \geq 2)$ with $[\alpha_{\B}]_{\F} \in \E_M$.

Now let $d(t) = (d_1(t_1), \ldots, d_g(t_g))^t$, where $d_j(t_j) = \pi_{\B} u_j t_j + t_j^p$.  Then $d(t) \in \E_M$ by construction, and hence, by Lemma \ref{Haz13.3.3}, there exists precisely one $g$-dimensional formal group $F_d$ over $\OL_{\B}$ such that $d$ is an endomorphism of $F_d$.  We will construct a formal group $\G$, which is a product of one-dimensional Lubin-Tate formal groups, with $d(t)$ an endomorphism of $\G$, and hence, by Lemma \ref{Haz13.3.3}, $\G = F_d$.  

Now apply Lemma \ref{Haz13.3.3} with $M = \pi_{\B} u_j$ for any $1 \leq j \leq g$.  
Then $d_j(t_j) \in \E_{\pi_{\B} u_j}$, so there exists precisely one one-dimensional Lubin-Tate type formal group $\G_j$ over $\OL_{\B}$ such that $d_j$ is an endomorphism of $\G_j$ for each $1 \leq j \leq g$.  Let $\G = \oplus_{j=1}^g \G_j$.  Therefore $\G$ is a $g$-dimensional commutative higher dimensional Lubin-Tate formal group.  Furthermore, $d(t) = (d_1(t_1), \ldots, d_g(t_g))^t$ is an endomorphism of $\G$ by construction, and hence $\G = F_d$.  By Lemma \ref{Haz13.3.3}, we have that $\F$ is strictly isomorphic to $\G$ over $\OL_{\B}$.

(ii)  For $1 \leq j \leq g$, we know that $\G_j$ is a one-dimensional Lubin-Tate formal group over $\OL_{\B}$ such that $d_j(t_j) = \phi_j(\alpha_{\B})t_j + t_j^p$ is an endomorphism of $\G_j$.  Following \cite{Lubin-Tate}, since $\phi_j(\alpha_{\B})$ is a uniformizer at $\B$, we define an endomorphism  $[\phi_j(\alpha_{\B})]_{\G_j}(t_j) := d_j(t_j)$, and for any $\gamma \in \OL_{\B}$, let  $[\gamma]_{\G_j}(t_j)$ be the unique power series such that $[\gamma]_{\G_j}(t_j) = \gamma t_j + (d^{\circ} \geq 2)$ and  $[\phi_j(\alpha_{\B})]_{\G_j} ([\gamma]_{\G_j}(t_j)) = [\gamma]_{\G_j}([\phi_j(\alpha_{\B})]_{\G_j}(t_j))$.  Let $m \in \Z \subset \OL_{\B}$, then it is easy to see that multiplication by $m$ coming from the group law on $\G_j$ is an endomorphism of $\G_j$, which must commute with any endomorphism, so is equivalent to the endomorphism $[m]_{\G_j}$ defined by Lubin-Tate theory.  Therefore we have an embedding of $\Z[\alpha_{\B}]$ into $\End(\G_j)$ given by $\beta \mapsto [\phi_j(\beta)]_{\G_j}$ for $\beta \in \Z[\alpha_{\B}]$, since $\phi_j(m) = m$ for $m \in \Z$.  Furthermore, the map $\beta \mapsto \oplus_{j=1}^g [\phi_j(\beta)]_{\G_j}$ gives an embedding of $\Z[\alpha_{\B}]$ into $\End(\G)$.  Since $\B$ is a prime of $K$ of degree one over $\Q$ by assumption, we have that $\Z[\alpha_{\B}]$ is of finite index in $\OK$ \cite[pg. 88]{LangCM}.  Therefore we can extend the above map to an embedding of $\OK$ into $\End(\G)$, where the jacobian of $[\alpha]_{\G} = \diag \langle \phi_1(u), \ldots, \phi_g(u) \rangle$.

Now, by (i), there is a strict isomorphism between $\F$ and $\G$, which we denote by $\beta$.
Therefore $\beta \circ [\alpha]_{\F} \circ \beta^{-1}$ is an endomorphism of $\G$ for all $\alpha \in \OK$.  Since $\beta$ is a strict isomorphism, $j(\beta) = j(\beta^{-1})$ is the identity matrix, and hence $j(\beta \circ [\alpha]_{\F} \circ \beta^{-1}) = j(\beta)j([\alpha]_{\F})j(\beta^{-1}) = j([\alpha]_{\F})$.  Thus the kernel of $[\alpha]_{\F}$ can be identified with the kernel of $\beta \circ [\alpha]_{\F} \circ \beta^{-1}$.  Since these jacobians are equal, Lemma \ref{Kpcmact} shows us that $\beta \circ [\alpha]_{\F} \circ \beta^{-1}(t) =  Nt + (d^{\circ} \geq 2)$, where $N = \diag \langle \phi_1(\alpha), \ldots, \phi_g(\alpha) \rangle$.  Therefore, for $\alpha \in \OK$, $[\alpha]_{\G}$ has the same jacobian as $\beta \circ [\alpha]_{\F} \circ \beta^{-1}$, and hence it is easy to see that they have the same kernel.  Thus we can identify $\F[\alpha]$ with $\G[\alpha]$ for any $\alpha \in \OK$.  

(iii)  By Corollary \ref{AFtors}, we know that we can identify $A[\p^n]$ with $\F[\pi^n]$.  By construction of the embedding of $\OK \hookrightarrow \End(\G)$ in (ii), we can identify $\F[\pi^n]$ with $\G[\pi^n]$.  By definition, we have
\begin{equation*}
\begin{split}
\G[\pi^n] & = \mbox{ker} \{ [\pi^n]_{\G}: \G(\M) \rightarrow \G(\M) \} \\
          & = \mbox{ker} \left \{ \oplus_{j = 1}^g [\phi_j(\pi^n)]_{\G_j}: \oplus_{j=1}^g \G_j(\M) \rightarrow \oplus_{j=1}^g \G_j(\M) \right \} \\
          & = \oplus_{j=1}^g \G_j[\phi_j(\pi^n)] \\
          & \cong \G_i[\phi_i(\pi^n)] = \G_i[\pi_i^n].
\end{split}
\end{equation*}
\noindent
Indeed, for $j \neq i$, $\ord_{\B} \phi_j(\pi^n) = 0$.  Therefore $[\phi_j(\pi^n)]_{\G_j}$ is an automorphism, and hence $\G_j[\phi_j(\pi^n)] = \{ O \}$.
\end{proof}

\begin{cor} \label{piiaction}
Let everything be as in Proposition \ref{fmlgrpdecomp}, then 
\begin{equation*}
[\pi_i]_{\G_i}(t_i) = \pi_i t_i + u t_i^p + \pi_i \alpha + \beta,
\end{equation*}
\noindent
where $u$ is unit in $\OL_{\B}$, $\alpha$ and $\beta$ are power series in $t_i$ with lowest terms of degree two and $2p$ respectively.  
\end{cor}

\begin{proof}
By Proposition \ref{fmlgrpdecomp}, we know that $[\alpha_{\B}]_{\G} = \oplus_{j=1}^g [\phi_i(\alpha_{\B})]_{\G_j}$.  By Lemma \ref{frobdecomp}, we know that $\ord_{\B} \phi_i(\alpha_{\B}) = 1$, and hence $\pi_i = \phi_i(\alpha_{\B}) u_i$ with $u_i$ a unit in $\OL_{\B}$.  Then we have the following (dropping the $\G_i$ in the notation).
\begin{equation*} 
\begin{split}
[\pi_i](t_i) & = [u_i]([\phi_i(\alpha_{\B})](t_i)) \\
             & = [u_i](\pi_i u_i^{-1} t_i + t_i^p) \\
             & = u_i (\pi_i u_i^{-1} t_i + t_i^p) + (d^{\circ} \geq 2)(\pi_i u_i^{-1} t_i + t_i^p) \\
             & = \pi_i t_i + u_i t_i^p + (d^{\circ} \geq 2)(\pi_i u_i^{-1} t_i + t_i^p).    
\end{split}
\end{equation*}
\end{proof}

\section{Division fields}\label{sec:divfld}

In this section, we are interested in algebraic properties of the field extensions of $K$ generated by adjoining $\p^n$-torsion points to $K$ for $K$ and $\p$-torsion from $A$ as in section \ref{sec:setup}.  (Recall that this means that $\p$ is a first-degree prime of $K$.)  We will denote \emph{the field of $\p^n$-division values of $K$} by $K_n = K(A[\p^n])$.

\begin{prop} \label{KnKstruct}
Let $\B$ be a prime of $K$ lying over $p$.
\begin{itemize}
\item[(i)] $K_n/K$ is totally ramified at $\B \in S$, and unramified at $\B \notin S$.  
\item[(ii)]  $K_n/K$ is a cyclic extension of degree $p^{n-1}(p-1)$.  
\item[(iii)]  Let $\B = \phi_i(\p)$ for fixed a fixed $i$ with $1 \leq i \leq g$ and let $\B_n$ be the unique prime of $K_n$ above $\B$.  Then $\ord_{\B_n} t_i (v) = 1$ for any $v \in A[\p^n] \setminus A[\p^{n-1}]$. 
\end{itemize}
\end{prop}

\begin{notation}
We denote by $S_n$ the collection of primes $\B_n$ of $K_n$ described in (iii).  By Proposition 5.3 of \cite{LangCM} and (ii),  we fix isomorphisms 
\begin{equation*}
\gal(K_n/K) \cong (\OK/\p^n \OK)^{\times} \cong (\Z / p^n \Z)^{\times}.
\end{equation*}
\end{notation}

\begin{proof}[Proof of Proposition \ref{KnKstruct}]
For now, fix a prime $\B = \phi_i(\p)$ with $1 \leq i \leq g$.  By Proposition \ref{fmlgrpdecomp}, we can identify $A[\p^n]$ with $\G_i[\phi_i(\pi^n)]$. 

On the one hand, if $\B_n$ is a prime of $K_n$ lying over $\B$, we can make the following identification:   
\[(K_n)_{\B_n} = (K(A[\p^n]))_{\B_n} = K_{\B}(A[\p^n]) = K_{\B}(\G_i[\phi_i(\pi^n)]).\]
\noindent
Let $F = K_{\B}$, and $F_n = K_{\B}(\G_i[\phi_i(\pi^n)])$.  Since $\G_i$ is a Lubin-Tate formal group, $F_n/F$ is a totally ramified cyclic extension of degree $p^{n-1}(p-1)$ (see \cite{Lubin-Tate}), and hence $[K_n:K] \geq p^{n-1}(p-1)$.  

On the other hand, it is well known that $\gal(K_n/K)$ is isomorphic to a subgroup of $(\OK/\p^n \OK)^{\times}$.  Since $p$ splits completely, we have that $(\OK/\p^n \OK)^{\times} \cong (\Z / p^n \Z)^{\times}$, where the latter is a cyclic group of order $p^{n-1}(p-1)$.  Therefore $[K_n:K] \leq p^{n-1}(p-1)$.  This completes the proof of half of (i) and all of (ii).  

For (iii), let $\lambda_n = t_i(v)$ for some $v \in A[\p^n] \setminus A[\p^{n-1}]$.  Since $\G_i$ is a Lubin-Tate formal group, $F_n = F(\lambda_n)$ and $N_{F_n/F}(-\lambda_n) = \phi_i(\pi)$ (see \cite[pg. 348]{Neukirch}).  Note that $(\phi_i(\pi)) = \B \OL_{\B}$, and therefore $\ord_{\B_n} \lambda_n = 1$.  

In order to finish (i), let $\B \notin S$ and $I_{\B_n}$ be the inertia group of $\B_n/\B$.  Assume that $\B$ ramifies in $K_n$, then $I_{\B_n} \neq 1$, i.e., there exists a $\sigma \neq 1 \in I_{\B_n}$ such that $\sigma(\alpha) \equiv \alpha \bmod \B_n$ for all $\alpha \in \OL_{K_n}$.  Therefore $\widetilde{\sigma(a)} = \widetilde{a}$ for $a \in A[\p^n] \setminus A[\p^{n-1}]$.  Since $\sigma \neq 1$, $\sigma(a) - a = b \in A[\p^n] \setminus \{ O \}$.  Recall that $A$ has good reduction at all primes above $p$, therefore the addition map reduces to a morphism on $\tilde{A}$ \cite[pg. 271]{Hindry-Silverman}.  Hence reduction commutes with addition, and we have that  
\[ \tilde{O} = \widetilde{\sigma(a)} - \widetilde{a} = \widetilde{\sigma(a) - a}  = \tilde{b}. \]
\noindent
This says that $b$ is in the kernel of reduction, but this contradicts Lemma \ref{primetop}.  Thus $I_{\B_n} = 1$, and $\B$ does not ramify as we claimed.
\end{proof}

Now we are in a position to determine some of the arithmetic properties of the fields $K_n/K$.

\begin{lem} \label{Kncond}
Let $\B_1 \in S_1$ and let $\chi$ be a non-trivial character on $\gal(K_n/K_1)$ with $n > 1$.  Then $\ord_{\B_1} \cond(\chi) \geq p$.
\end{lem}
   
\begin{proof}
Let $\chi$ be a non-trivial character on the group $\gal(K_n/K_1)$, which is a cyclic group of order $p^{n-1}$ by Proposition \ref{KnKstruct}.  For ease of notation, let $L_{\chi}$ denote the fixed field of $\chi$.  Then $L_{\chi} = K_r$ for some $2 \leq r \leq n$, and hence $\cond(\chi) = \cond(K_r/K_1)$ and $\cond(K_2/K_1)|\cond(K_r/K_1)$.  Therefore it is enough to show that $\cond((K_2)_{\B_2}/(K_1)_{\B_1}) = \B_1^p$.   Once we calculate $D((K_2)_{\B_2}/(K_1)_{\B_1})$, the result follows by applying the conductor-discriminant formula (see \cite[pg. 534]{Neukirch}).  

Let $v \in A[\p^2] \setminus A[\p]$, then by Proposition \ref{KnKstruct}, $t_i(v)$ is a uniformizer at $\B_2$, and hence $\OL_{(K_2)_{\B_2}} = \OL_{(K_1)_{\B_1}}[t_i(v)]$.  Let $f$ be a minimum polynomial for $t_i(v)$, then
\begin{equation} \label{discminpolyform}
\begin{split}
D((K_2)_{\B_2}/(K_1)_{\B_1}) 
  & = N_{(K_2)_{\B_2}/(K_1)_{\B_1}} (f'(t_i(v))  \\
  & = N_{(K_2)_{\B_2}/(K_1)_{\B_1}}  (\prod_{\stackrel{\scriptstyle{u = v + w}}{w \in A[\p] \setminus \{ O \}}} (t_i(v) - t_i(u)) )  \\
  & = \prod_{\stackrel{\scriptstyle{u, v \in A[\p^2] \setminus A[\p]}}{u \neq v}} (t_i(v) - t_i(u)) \\
  & = \B_1^{p(p-1)}.
\end{split}
\end{equation}
\noindent
Indeed, for each $v \neq u \in A[\p^2] \setminus A[\p]$, there exists $w \in A[\p] \setminus \{ O \}$ such that $t_i(v) - t_i(u) = t_i(u + w) - t_i(u) =  t_i(w) + (d^{\circ} \geq 2)(t_i(u), t_i(w))$.  By Proposition \ref{KnKstruct}, $\ord_{\B_1} t_i(w) = 1$, and hence by comparing terms of least valuation we have $\ord_{\B_1} (t_i(v) - t_i(u)) = 1$.  

Since $(K_2)_{\B_2}/(K_1)_{\B_1}$ is an extension of degree $p$, the conductor-discriminant formula gives us that $D((K_2)_{\B_2}/(K_1)_{\B_1}) = \cond((K_2)_{\B_2}/(K_1)_{\B_1})^{p-1}$, and comparing this with $\eqref{discminpolyform}$ completes the result.
\end{proof}

\begin{prop} \label{gdred}
$A$ has everywhere good reduction over $K_1$.  
\end{prop}

\begin{proof}
Since $A$ has complex multiplication by $\OK$, good reduction at all primes lying over $p$, and $p$ splits completely in $K$, the result follows from slight modification to the proof of Theorem 2 of \cite{Coates-Wiles} (see \cite{Rowe} for more details).
\end{proof}

\section{The tower}\label{sec:tower}

Let $h$ denote the class number of $K$.  Assume that the $\OK$-rank of $A(K)$ is positive.  Let $Q$ be a point of infinite order of $A$ rational over $K$ such that $Q \in A(K) \setminus \p A(K)$, where $\p A(K) = \cup_{\alpha \in \p} [\alpha]A(K)$ (this set is non-empty since the free part of $A(K)$ is a finitely generated $\OK$-module).  

Let $m = rh$ be a positive integer multiple of the class number of $K$.  
We have that $\p^{m} = (\gamma)$ for some $\gamma \in \OK$ (hence $A[\p] = A[\gamma]$).  Let $Q_m \in A(\overline{K})$ be such that $[\gamma]Q_m = Q$.  We define $L_m = K_m(Q_m)$ to be the field of \emph{$\p$-division values of $Q$} which is independent of the choice of $Q_m$ and $\gamma$.  Then, for $1 \leq n \leq m$, we define $L_n = K_n( \cup_{\alpha \in \p^{m-n}} [\alpha]Q_m)$ to be the field of \emph{$\p^n$-division values of $Q$}.

\begin{rems} ~
\begin{enumerate}
\item Since the class number of $K$ is often not one, we cannot define the field of $\p^n$-division values of $Q$ in the usual way.  Therefore we use the fact that $\OK$ is a Dedekind domain to note that $\p^{m-n}$ is generated by two elements of $\OK$, say $\p^{m-n} = (\alpha_1, \alpha_2)$.  Then it is not hard to see that $L_n = K_n([\alpha_1]Q_m, [\alpha_2]Q_m)$ regardless of the choice of generators for $\p^{m-n}$.

\item This definition is independent of our choice of $m$.  Let $m_1 = r_1 h$ and $m_2 = r_2 h$ with $n < m_1 < m_2$.  We can choose $Q_{m_2}$ recursively so that $[\gamma]Q_{r_2 h} = Q_{(r_2 -1)h}$, and hence $[\gamma^{r_2-r_1}]Q_{m_2} = Q_{m_1}$.  Then $\p^{m_2 -n} = \p^{r_1 h -n} \p^{(r_2 - r_1)h} = \p^{m_1 - n} (\gamma^{r_2 - r_1})$ and  
\begin{equation*}
K_n(\cup_{\alpha \in \p^{m_1 -n}} [\alpha]Q_{m_1}) 
 =  K_n ( \cup_{\alpha \in \p^{m_1 -n}} [\alpha][\gamma^{r_2-r_1}]Q_{m_2}) 
 =  K_n ( \cup_{\beta \in \p^{m_2 - n}} [\beta]Q_{m_2}).
\end{equation*}
\item From now on, we may assume that $m$ is a suitably large, fixed multiple of $h$.   Via our isomorphism $\OK/\p^n \OK \cong \Z / p^n \Z$, we indentify $A[\p^n] \subset A[\p]$ with the subgroup $\p^{m-n} \OK/ \p \OK \subset \OK/\p \OK$.  
\end{enumerate}
\end{rems}

Following \cite{Gupta}, we call the fields $L_n \supset K_n \supset K$ the Coates-Wiles tower.  It is easy to see that $L_n/K$ is Galois for all $n$.  Furthermore, $G_m = \gal(L_m/K)$ is  identifiable with a subgroup of 
\begin{equation*}
H_m = \left \{ \begin{pmatrix}1 & * \\ 0 & * \end{pmatrix} \in GL_2(\OK/\p \OK) \right \}.
\end{equation*}
\noindent
The action of $\displaystyle \sigma = \begin{pmatrix} 1 & a_{\sigma} \\ 0 & b_{\sigma} \end{pmatrix} \in G_m$ is given by
\begin{equation*}
\begin{split}
\sigma(Q_m) & = Q_m + a_{\sigma}, \\
\sigma(R)   & = [b_{\sigma}]R,
\end{split}
\end{equation*}
\noindent
where $R, a_{\sigma} \in A[\p] \cong \OK / \p$ and $b_{\sigma} \in (\OK/\p)^{\times}$.

Since $\p$ splits completely and we have complex multiplication, we expect the extension $L_n/K$ to be ``as big as possible''.  

\begin{lem} \label{LmKmIso}
$\gal(L_m/K_m) \cong A[\p]$.
\end{lem}

\begin{proof}
We will first show that $L_1 \neq K_1$.  Let $\gamma$ and $\pi$ be as per usual, so that $\p = (\gamma)$ and $\ord_{\p} \pi = \ord_p \pi = 1$.  Then $\p^{m-1} = (\gamma, \pi^{m-1})$, and we have that $L_1 = K_1([\gamma]Q_m, [\pi^{m-1}]Q_m) = K_1([\pi^{m-1}]Q_m)$.  For $\tau \in \gal(L_1/K)$, we have $\tau([\pi^{m-1}]Q_m) = [\pi^{m-1}]Q_m + a_{\tau}$ with $a_{\tau} \in A[\p] \setminus \{ O \}$.  

Now assume that $L_1 = K_1$, and define the map $f(\tau) = a_{\tau}$.  It is clear that this defines a 1-cocycle of $\gal(L_1/K) = \gal(K_1/K)$ into $A[\p]$.
Since $A[\p]$ is a $p$-group, $H^1 (\gal(K_1/K), A[\p])$ must be a $p$-group.  But Proposition \ref{KnKstruct} shows that $[K_1:K] = p-1$  and \cite[pg. 130]{SerreLocFld} shows that $H^1 (\gal(K_1/K), A[\p])$ is annihilated by multiplication by $p-1$.  Therefore $H^1(\gal(K_1/K), A[\p]) = 0$ and $f$ must be a coboundary.  Thus there exists $c \in A[\p]$ such that $f(\tau) = \tau c - c$ for all $\tau \in \gal(K_1/K)$.  Therefore $\tau$ fixes $[\pi^{m-1}]Q_m - c$ for all $\tau \in \gal(K_1/K)$, and hence $[\pi^{m-1}]Q_m - c \in A(K)$.  By our choice of $\pi$, we must have $\pi = \delta \gamma$ with $\delta$ and $p$ relatively prime.  But then $[\pi]([\pi^{m-1}]Q_m -c) = [\pi]Q_m = [\delta \gamma]Q_m = [\delta]Q \in A(K)$, which shows that $Q \in [\p]A(K)$.  This is a contradiction of the choice of $Q$, and hence $L_1 \neq K_1$.

Let $\sigma \in \gal(L_m/K_m)$.   Then it is not hard to see that $\sigma Q_m = Q_m + a_{\sigma}$ with $a_{\sigma} \in A[\p]$ and the map $\sigma \mapsto a_{\sigma}$ gives us an embedding of $\gal(L_m/K_m)$ into $A[\p]$ (for more details see \cite{Rowe}).  Thus the map $\sigma \mapsto a_{\sigma}$ injects $\gal(L_1/K_1)$ into $A[\p]$.  Since $L_1/K_1$ is non-trivial, $\gal(L_1/K_1)$ is isomorphic to a non-trivial subgroup of $A[\p]$.  But the only non-trivial subgroup of $A[\p]$ is $A[\p]$.  

Since $L_1/K_1 \cong A[\p]$, the rest lemma follows by modifying an argument of Lang for elliptic curves \cite{Lang2}.     
\end{proof}

Now we show that, for $n=m$, $L_m/K$ is as big as possible.

\begin{prop} \label{LmKrep}
$G_m = H_m$.
\end{prop}

\begin{proof}
By Proposition \ref{KnKstruct}, we have that $[K_m:K] = p^{m-1}(p-1)$.  By Proposition 5.3 of \cite{LangCM}, we have that $A[\p] \cong \OK/\p \OK$, and hence $\#(A[\p]) = p$.  By Galois theory and Lemma \ref{LmKmIso}, we must have $\#(G_m) = [L_m:K] = p p^{m-1} (p-1)$.  By construction $G_m \hookrightarrow H_m$, and a simple calculation shows that $\#(H_m) = p p^{m-1}(p-1) = \#(G_m)$.  Therefore $G_m = H_m$.
\end{proof}

\begin{notation} 
Recall that we fixed an isomorphism such that $A[\p] \cong \OK/\p \OK \cong \Z / p \Z$, where $A[\p^{m-n}] \subset A[\p]$ corresponds to $\p^n \OK/ \p \OK$.  For $\alpha \in \OK/\p\OK$, we say $\alpha \equiv 0 \bmod \p^n$ if $\alpha$ comes from $A[\p^{m-n}]$.  We will say that $\beta \equiv 1 \bmod \p^n$ if $\beta \in \left( \OK/\p \OK \right)^{\times}$ acts on $A[\p^n] \subset A[\p]$ as the identity automorphism.    
\end{notation}

Since Proposition \ref{LmKrep} gives us the matrix representation of $\gal(L_m/K)$, we can use properties of this matrix group to try and determine the representation of its subgroups.

\begin{prop} \label{LnKnrep}
$\gal(L_n/K_n) \cong A[\p^n]$, and hence $L_n/K$ is ``as big as possible''.
\end{prop}

\begin{proof}
It is not hard to see that $L_n \cap K_m = K_n$, and that we have the following isomorphisms:
\begin{equation} \label{LmKnrep}
\gal(L_m/K_n) \cong \left \{ \left( \begin{smallmatrix} 1 & \ast \\ 0 & \beta \end{smallmatrix} \right) \in \mbox{GL}_2 \left( \OK/\p \OK \right) : \beta \equiv 1 \bmod \p^n \right \}
\end{equation}
and
\begin{equation} \label{LmLnrep}
\gal(L_m/L_n) \cong \left \{ \left( \begin{smallmatrix} 1 & \alpha \\ 0 & \beta \end{smallmatrix} \right) \in \mbox{GL}_2 \left( \OK/\p \OK \right) : \alpha \equiv 0, \beta \equiv 1 \bmod \p^n \right \}.
\end{equation}
\noindent
From equations $\eqref{LmKnrep}$ and $\eqref{LmLnrep}$, we can compute that $[L_m:K_n] = p p^{m-n}$ and $[L_m:L_n] = p^{2m -2n}$.  Therefore $[L_n:K_n] = p^{2m -n}/p^{2m-2n} = p^n$.

Since $A[\p^n]$ is a cyclic group of order $p^n$ and $[L_n:K_n] = p^n$, we need only show that $\gal(L_n/K_n)$ is isomorphic to a cyclic group.  Since $L_n \cap K_m = K_n$, we have $\gal(L_n/K_n) \cong \gal(L_n K_m/K_m)$.  But the second group is a quotient of the cyclic group $\gal(L_m/K_m)$, hence it is cyclic.  
\end{proof}

\begin{figure}[htbp]
\begin{center}
\epsfig{file=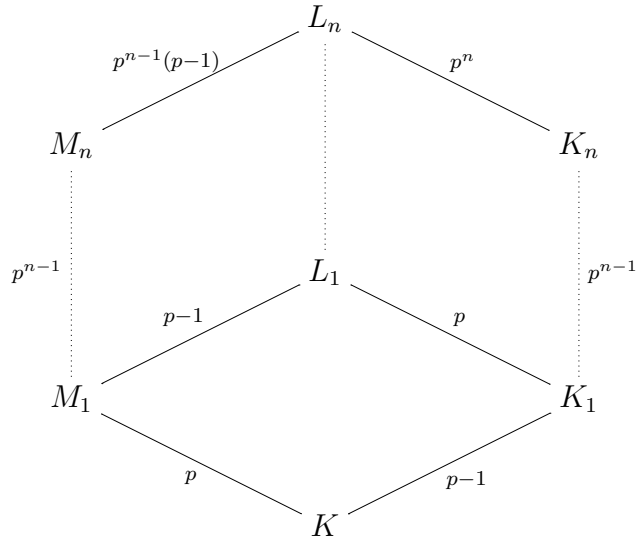}
\caption{The Coates-Wiles tower}
\label{fig:cwtower}
\end{center}
\end{figure}

In order to compute the conductor of $L_n/K_n$, we will need to compare calculations made up different branches of the Coates-Wiles tower (see Figure \ref{fig:cwtower}).  Moreover, we will find the following subextensions of $L_n/K$ helpful in our calculations: 
\begin{equation}\label{def:Mn}
M_n := K(\cup_{\alpha \in \p^{m-n}} [\alpha]Q_m)
\end{equation}
and
\begin{equation}\label{def:Mntilde}  
\widetilde{M}_n := K_1(\cup_{\alpha \in \p^{m-n}} [\alpha]Q_m) = K_1M_n \mbox{ for $1 \leq n < m$.}
\end{equation}

\section{The Filtration.} \label{sec:filt}

We may assume that $Q \in A^{\circ}(K_{\B})$ for all $\B \in S$, since after multiplying $Q$ by a suitable integer relatively prime to $p$ this is the case.  Let $\pi \in \OK$ such that $\ord_p \pi = \ord_{\p} \pi = 1$, and recall that $\p = (\gamma)$.  Then $\pi = \gamma \delta$ with $\delta$ and $p$ relatively prime and $\p^{m-n} = (\gamma, \pi^{m-n})$ for $1 \leq n \leq m$.  Thus we have $L_n = K_n([\pi^{m-n}]Q_m)$.  

We are interested in computing $\cond(L_n/K_n)$, and the following lemma shows that $\cond(L_n/K_n)$ is only divisible by primes of $S_n$.

\begin{lem} \label{LmKmunram}
$L_n/K_n$ is unramified outside of $S_n$. 
\end{lem}

\begin{proof}
In Lemma \ref{gdred}, we showed that $A$ has everywhere good reduction over $K_1$, and hence over $K_n$.  It follows easily from Lemma \ref{primetop} that $L_n$ can only be ramified at primes of $K_n$ which lie above $p$ (see \cite{Serre-Tate}).  

First we will show the result for $n=m$.  Let $\q \notin S_m$ be a prime of $K_m$ above $p$, $F = (K_m)_{\q}$ the completion of $K_m$ at $\q$, and    $A^{\circ}(F)$ the kernel of reduction $\bmod$ $\q$.  

Let $\bigq$ be a prime of $\overline{K_m}$ lying over $\q$, and $I_{\bigq}$ the inertia group of $\bigq$ over $\q$.  Then $I_{\bigq}$ acts trivially on $Q_m$ if and only if $L_m$ is unramified at $\q$.  Assume $I_{\bigq}$ does not act trivially, then there exists $\sigma \neq 1 \in I_{\bigq}$ such that $\sigma(Q_m)$ and $Q_m$ reduce to the same element $\bmod$ $\q$.  Therefore $\sigma(Q_m) - Q_m \in A^{\circ}(F)$, but $\sigma(Q_m) - Q_m = \alpha \in A[\p]$.  Indeed,%
\begin{equation*}
[\gamma](\sigma(Q_m) - Q_m) = \sigma([\gamma]Q_m) - [\gamma]Q_m = \sigma(Q) - Q = O,
\end{equation*}
\noindent
where $\p = (\gamma)$ with $\gamma \in \OK$.  By good reduction at $\q$, we know that reduction commutes with addition, and hence we have
\begin{equation*}
\widetilde{O} = \widetilde{\sigma(Q_m) - Q_m} = \widetilde{\sigma(Q_m)} - \widetilde{Q}_m = \widetilde{\alpha}.
\end{equation*}
\noindent
By Lemma \ref{primetop}, $\alpha = O$, and hence $\sigma$ acts trivially on $Q_m$.  This contradiction completes the case of $n=m$.  For $1 \leq n < m$, let $\q' \notin S$ be a prime of $K$ lying over $p$.  Then the case of $n=m$ and Proposition \ref{KnKstruct} combine to show us that $\q$ does not ramify in $L_m$.  Therefore $\q$ does not ramify in $L_n \subset L_m$, and hence $L_n/K_n$ is unramified outside of $S_n$.    
\end{proof}

For what follows, it is important to recall that  
\begin{enumerate}
\item $\B$ is a prime of $K$ such that $\B = \phi_i(\p)$ with $1 \leq i \leq g$, where $\phi_i$ refers to an element of the CM type of $A$, and 

\item $(\phi_i(\pi)) = (p) = \B \OL_{\B}$ locally for some $\pi \in \OK$.
\end{enumerate}

By Proposition \ref{fmlgrpdecomp}, we are able to define a filtration of $A^{\circ}(K_{\B})$ for $\B \in S$.  Let $A^{\circ}(K_{\B}) = A_0(K_{\B})$ and define
\begin{equation} \label{filtdef}
A_j(K_{\B}) = \{ R \in A^{\circ}(K_{\B}) : \ord_{\B} t_i(R) > j \};
\end{equation}
\noindent
so $A_0(K_{\B}) \supseteq A_1(K_{\B}) \supseteq A_2(K_{\B}) \supseteq \cdots$

\begin{prop} \label{filtration}
Let $L$ be a finite extension of $K$ with prime $P$ lying over $\B \in S$ and let $\beta \in \OK$ be relatively prime to $p$.  Then 

\begin{itemize}
\item[(i)]  $[\pi]_{\G}: A_j(K_{\B}) \rightarrow A_{j+1}(K_{\B})$ is an isomorphism for all $j \geq 0$.  

\item[(ii)]  $[\beta]_{\G}: A^{\circ}(L_{P}) \rightarrow A^{\circ}(L_P)$ is an isomorphism.     
\end{itemize}
\end{prop}

\begin{proof}
(i)  Proposition \ref{fmlgrpdecomp} shows that $\G = \oplus_{k=1}^g\G_k$, where $\G_k$ is a one-dimensional Lubin-Tate formal group.  What we have done is to construct a filtration on the $i$th component of this formal group, and by \cite{SerreLieAlg} we have an isomorphism on the $i$th component.  Now since $[\pi]_{\G} = \oplus_{k=1}^g[\phi_k(\pi)]_{\G_k}$ and $\ord_{\B} \phi_k(\pi) = 0$ for $k \neq i$, $[\phi_k(\pi)]_{\G_k}$ is an automorphism of $\G_k(\B)$.  Indeed, $j([\phi_k(\pi)])$ is a unit in $\OL_{\B}$ for $k \neq i$.  Since $A_0(K_{\B}) \cong \G(\B) = \oplus_{k=1}^g \G_i(\B)$,  $[\pi]_{\G}$ is an isomorphism of $A_j(K_{\B})$ onto $A_{j+1}(K_{\B})$.

(ii)  Given $\beta \in \OK$ with $\beta$ and $p$ relatively prime, we have $[\beta]_{\G} = \oplus_{k=1}^g[\phi_i(\beta)]_{\G_k}$, and $\phi_k(\beta)$ is a unit in $K_{\B}$ for all $1 \leq k \leq g$.  Then the result follows from a slight modification to the proof of Theorem 3 \cite[LG 4.25]{SerreLieAlg}. 
\end{proof}

\begin{notation} 
Let $e_{\B}$ be the smallest integer such that $Q \notin A_{e_{\B}}(K_{\B})$, $e = \mbox{min}_{\B}\{e_{\B}\}$, and $\SOP = \{ \B \in S : e_{\B} = e \}$.
\end{notation}

Let $P_n$ be a prime of $L_n$ lying over $\B_n$, for fixed $\B \in S$.  Let $Q = \widetilde{Q}_0$ and choose $\widetilde{Q}_n$ recursively such that $[\pi]\widetilde{Q}_n = \widetilde{Q}_{n-1}$.  


\begin{prop} \label{Lnlocal}
Let $P_n$ be a prime of $L_n$ lying over $\B_n \in S_n$, then
\begin{equation*}
(L_n)_{P_n} = (K_n)_{\B_n}(t_i(\widetilde{Q}_n)).
\end{equation*}
\end{prop}  

\begin{proof}
We will first show that $(L_m)_{P_m}= (K_m)_{\B_m}(\widetilde{Q}_m)$.  By definition,   
we have that $(L_m)_{P_m} = (K_m)_{\B_m}(Q_m)$.  Since $Q \in A^{\circ}(K_{\B})$ and $[\gamma]Q_m = Q$, there is some choice of $Q_m \in A^{\circ}((L_m)_{P_m})$.  Since $A[\p]$ is in the kernel of reduction, we actually have all choices of $Q_m \in A^{\circ}((L_m)_{P_m})$.  By $\eqref{kr}$, we have that $A^{\circ}((L_m)_{P_m}) \cong \G(P_m)$, where the map is given by $R \mapsto t(R) = (t_1(R), \ldots, t_g(R))^t$.  Hence
\begin{equation*}
(K_m)_{\B_m}(Q_m) = (K_m)_{\B_m}(t_1(Q_m), \ldots, t_g(Q_m))
\end{equation*}
\noindent  
Note that $[\pi]\widetilde{Q}_m = [\gamma] ([\delta]\widetilde{Q}_m) = Q$, and hence $[\delta]\widetilde{Q}_m = Q_m + a$ with $a \in A[\p]$.  Therefore, in the formal group, we have that $[\delta]_{\G}t(\widetilde{Q}_m) = t(Q_m + a) = t(Q_m) +_{\G} t(a)$.  Since $t(a) \in (K_m)_{P_m}$, we may assume without loss of generality that $a = O$.  Since $\delta$ is relatively prime to $p$, Proposition \ref{filtration} gives us that $[\delta]_{\G}$ is an automorphism , and hence invertible.  Therefore $[\delta]_{\G}t(\widetilde{Q}_m) = t(Q)$, and we have that
\begin{equation*}
(K_m)_{\B_m}(t_1(Q_m), \ldots, t_g(Q_m)) \subset (K_m)_{\B_m}(t_1(\widetilde{Q}_m), \ldots, t_g(\widetilde{Q}_m)).
\end{equation*}
\noindent
The reverse inclusion then comes from $[\delta^{-1}]_{\G}t(Q_m) = t(\widetilde{Q}_m)$.  For $1 \leq n < m$, recall that $L_n = K_n([\pi^{m-n}]Q_m)$.  As above,  we have that 
\begin{equation*}
(L_n)_{P_n} = (K_n)_{\B_n}([\pi^{m-n}]Q_m) = (K_n)_{\B_n}([\pi^{m-n}]\widetilde{Q}_m) = (K_n)_{\B_n}(\widetilde{Q}_n).
\end{equation*}
\noindent
Finally, recall that $[\pi^n]_{\G} = \oplus_{k=1}^g [\phi_k(\pi^n)]_{\G_k}$, and $[\phi_k(\pi^n)]_{\G_k}$ is an automorphism on $\G_k(\B)$ for $k \neq i$.  Hence $t_k(\widetilde{Q}_n) \in \G_k(\B)$ for $k \neq i$, which shows us that 
\begin{equation*}
(K_n)_{\B_n}(\widetilde{Q}_n) = (K_n)_{\B_n}(t_1(\widetilde{Q}_n), \ldots, t_g(\widetilde{Q}_n)) = (K_n)_{\B_n}(t_i(\widetilde{Q}_n)).
\end{equation*}
\end{proof}

\begin{notation}
Since $(K_n)_{\B_n}(\widetilde{Q}_n) = (K_n)_{\B_n}(t_i(\widetilde{Q}_n))$, we will use the notations interchangeably in what follows.  
\end{notation}

\begin{prop} \label{filt2}
Let $\B \in S$ and $e_{\B}$ be as above.   Then

\begin{itemize}
\item[(i)]  for $1 \leq n < e_{\B}$, $\B_n$ splits completely in $L_n$,

\item[(ii)]  $L_{e_{\B}}/K_{e_{\B}}$ is ramified over $\B_{e_{\B}}$, and

\item[(iii)]  $\B_{e_{\B}}$ splits completely in $L_{e_{\B}-1}K_{e_{\B}}$.  Hence $L_{e_{\B}}/ L_{e_{\B}-1}K_{e_{\B}}$ is a totally ramified extension of degree $p$ over any prime above $\B$.
\end{itemize}
\end{prop}

\begin{proof}
(i)  By the filtration $\eqref{filtdef}$, we see that there is a $Q_n^{'} \in A_{e_{\B}-n-1}(K_{\B})$ such that $[\pi^n]_{\G}t(Q_n^{'}) = t(Q)$.  But $[\pi^n]_{\G}t(\widetilde{Q}_n) = t(Q)$, hence $\widetilde{Q}_n = Q_n^{'} + u$ for some $u \in A[\p^n]$.  Since $Q_n^{'} + u \in A((K_n)_{\B_n})$, we see that $\B_n$ splits completely.  Indeed, there exists a prime $P_n$ of $L_n$ lying over $\B_n$ such that 
\begin{equation*}
(L_n)_{P_n}= 
(K_n)_{\B_n}(\widetilde{Q}_n) = (K_n)_{\B_n}(Q_n^{'} + u) = (K_n)_{\B_n}(u) = (K_n)_{\B_n}.
\end{equation*}  

(ii)  Since $p$ splits completely and we have the formal group decomposition of Proposition \ref{fmlgrpdecomp}, this follows from a slight modification to Theorem 11 of \cite{Coates-Wiles}.
  
(iii)  Let $n = e_{\B}$.  By (i), the decomposition group of $\gal(L_{n-1}/K_{n-1})$ for $\B_{n-1}$ is trivial.  It is not hard to see that $\gal(L_{n-1}/K_{n-1}) \cong \gal(L_{n-1}K_n/K_n)$.  Then the decomposition group of $\gal(L_{n-1}K_n/K_n)$ for $\B_n$ is isomorphic to a trivial group.  Indeed, these decomposition groups are the Galois groups of the corresponding local extensions.  Hence $\B_n$ splits completely in $L_{n-1}K_n$.
\end{proof}

We can reduce the computation of $\cond(L_n/K_n)$ to computing an Artin conductor.  
If $\chi$ is a character of a subgroup $H$ of $\gal(L_n/K)$, let $\chi^{\ast}$ be the corresponding induced character on $\gal(L_n/K)$, and $\char(n) = \{ \chi \in \widehat{\gal}(L_n/K_n) : \chi^{p^{n-1}} \neq 1 \}$, where $\widehat{\gal}(L_n/K_n)$ denotes the character group of $\gal(L_n/K_n)$.
The next lemma demonstrates the importance of computing $\cond(\chi)$ for $\chi \in \char(n)$.  

\begin{lem} \label{chare}
For any $\chi \in \char(e)$, $\cond(L_e/K_e) = \cond(\chi)$.
\end{lem}

\begin{proof}
By definition, $\chi^{p^{e-1}} \neq 1$ for $\chi \in \char(e)$.  Hence the fixed field of $\chi$, $L_{\chi}$, is not contained in $K_eL_{e-1}$, which is an extension of $K_e$ of degree $p^{e-1}$.  Therefore $L_{\chi} = L_e$.  Since $\cond(\chi) = \cond(L_{\chi}/K_e)$, we are done.
\end{proof}

\section{Conductor Calculations}\label{sec:cond}

Let $\SOP$ be as in the last section and for the rest of this section fix $\B \in \SOP$.  Further, we recall that $\B_n$ is the unique prime of $K_n$ sitting above $\B$ in $K$.  By Proposition \ref{filtration} and $\eqref{def:Mn}$, there exists a prime $\pee$ of $M_{e-1}$ over $\B$ such that 
\begin{equation*}
(M_{e-1})_{\pee} = K_{\B}.
\end{equation*}

Let $\pee_e$ be any prime of $L_{e-1}K_e$ lying over $\pee$ and $\pee_1$ its restriction to $\widetilde{M}_{e-1}$ (see $\eqref{def:Mntilde}$ to recall the definition of $\widetilde{M}_{e-1}$).  By $\eqref{filtdef}$ and Proposition \ref{filt2}, we see that both
\begin{equation*}
(L_{e-1}K_e)_{\pee_e} = (K_e)_{\B_e} \mbox{ and } (\widetilde{M}_{e-1})_{\pee_1} = (K_1)_{\B_1}.
\end{equation*}

\begin{figure}[htbp] 
\begin{center}
\epsfig{file=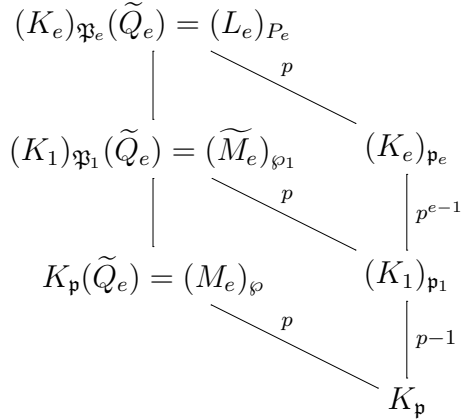}
\caption{The local tower}
\label{fig:localtower}
\end{center}
\end{figure}

Again by Proposition \ref{filt2}, there exists a unique totally ramified prime $P_e$ of $L_e$ over $\pee_e$ such that $[(L_e)_{P_e}:(L_{e-1}K_e)_{\pee_e}]=p$.  Let $\wp_1, \wp$ be the restrictions of $P_e$ to $\widetilde{M}_e$ and $M_e$ respectively.  The portion of the Coates-Wiles tower we are interested locally is shown in Figure \ref{fig:localtower}, where it is important to note that all of the extensions are totally ramified.

As a first step we will compute the conductor of $(K_1)_{\B_1}(\widetilde{Q}_e)/(K_1)_{\B_1}$, and then use this to compute the conductor of $(K_e)_{\B_e}(\widetilde{Q}_e)/(K_e)_{\B_e}$.  By Proposition \ref{filtration}, we have that $\ord_{\B}t_i(\widetilde{Q}_{e-1}) = 1$ and $t_i(\widetilde{Q}_{e-1}) = [\pi_i]_{\G_i} t_i(\widetilde{Q}_e)$.  In Corollary \ref{piiaction}, we saw that 
\begin{equation*}
[\pi_i]_{\G_i}(t_i) = \pi_i t_i + u t_i^p + \pi_i \alpha + \beta,
\end{equation*}
\noindent
where $\alpha$, $\beta$ are power series in $t_i$ with terms of lowest degree two and $2p$ respectively.  So
\begin{equation}\label{val1}
t_i(\widetilde{Q}_{e-1}) = [\pi_i]_{\G_i}(t_i)(\widetilde{Q}_e) = \pi_i t_i(\widetilde{Q}_e) + u t_i^p(\widetilde{Q}_e) + \pi_i\alpha(\widetilde{Q}_e) + \beta(\widetilde{Q}_e).
\end{equation}

Since $(M_e)_{\wp}/K_{\B}$ is a totally ramified extension of degree $p$, $\ord_{\wp} t_i(\widetilde{Q}_{e-1})=p$.  Note that $\ord_{\wp} \pi_i t_i(\widetilde{Q}_e) > p$, hence the term on the right-hand side of $\eqref{val1}$ with least valuation at $\wp$ can only be $u t_i(\widetilde{Q}_e)^p$.  Since $u$ is a unit in $\OL_{\B}$, we see that $\ord_{\wp} t_i(\widetilde{Q}_e) = 1$, and we also have $\ord_{\wp_1}(\widetilde{Q}_e) = p-1$ (since  $(K_1)_{\B_1}(\widetilde{Q}_e)/ K_{\B}(\widetilde{Q}_e)$ is a totally ramified extension of degree $p-1$).  Thus $t_i(\widetilde{Q}_e)$ has a $\wp_1$-adic expansion, $\sum_{i \geq p-1} a_i \eta_1^i$, with $a_i \in \{0, \ldots, p-1 \}$, $a_{p-1} \neq 0$, and $\eta_1$ a uniformizer at $\wp_1$.  Indeed, $\OL_{\wp_1}/\wp_1 \OL_{\wp_1} \cong \Z / p \Z$, so we can take $a_i \in \{0, \ldots, p-1 \}$.  Therefore we have that 
\begin{equation} \label{etauniform} 
(K_1)_{\B_1}(\widetilde{Q}_e)/(K_1)_{\B_1} = (K_1)_{\B_1}(\eta_1)/ (K_1)_{\B_1}.  
\end{equation}
\noindent
Let $\sigma \neq \tau \in \gal((K_1)_{\B_1}(\widetilde{Q}_e)/(K_1)_{\B_1})$, then $\sigma t_i(\widetilde{Q}_e)= t_i(\widetilde{Q}_e + u)$, $\tau t_i(\widetilde{Q}_e)= t_i(\widetilde{Q}_e + v)$ for some $u \neq v \in A[\p]$.  Using the formal group $\G_i$, we see
\begin{equation}\label{val2}
t_i(\widetilde{Q}_e + u) - t_i(\widetilde{Q}_e + v) = t_i(u-v) + (d^{\circ} \geq 2)(t_i(\widetilde{Q}_e + v), t_i(u-v)).
\end{equation}
\noindent 
By Proposition \ref{KnKstruct}, $t_i(u-v)$ is a uniformizer at $\B_1$, and hence  $\ord_{\wp_1} t_i(u-v) = p$.  We also have that $\ord_{\wp_1} t_i(\widetilde{Q}_e+v) = p-1$, because $\ord_{\wp} t_i(\widetilde{Q}_e+v) =\ord_{\wp}t_i(\widetilde{Q}_e) = 1$ (since any choice of $\widetilde{Q}_e$ will work in $\eqref{val1}$).  So by comparing the terms of $\eqref{val2}$ of least valuation at $\wp_1$ we see that 
\begin{equation}\label{val3}
\begin{split}
p & = \ord_{\wp_1} (t_i(\widetilde{Q}_e + u) - t_i(\widetilde{Q}_e + v)) \\
  & = \ord_{\wp_1}  (\sum_{i \geq p-1} a_i((\sigma \eta_1)^i - (\tau \eta_1)^i))  \\
  & = \ord_{\wp_1}  ((\sigma \eta_1 - \tau \eta_1)  \sum_{i \geq p-1} a_i  \sum_{j=0}^{i-1} (\sigma \eta_1)^j (\tau \eta_1)^{i-1-j})   \\
  & = \ord_{\wp_1} (\sigma \eta_1 - \tau \eta_1) + \ord_{\wp_1} ( \sum_{i \geq p-1} a_i \sum_{j=0}^{i-1} (\sigma \eta_1)^j (\tau \eta_1)^{i-1-j}).
\end{split}
\end{equation}
\noindent 
This follows because $x^n - y^n = (x-y)(\sum_{j=0}^{n-1} x^j y^{n-1-j})$.  Note that the second term on the last line on the right-hand side of $\eqref{val3}$ has $\ord_{\wp_1} \geq p-2$.  Indeed, $p-2$ is the smallest this valuation can be, since we get at least $p-2$ copies of terms with the same valuation as $\eta_1$ in the sum.  So $\eqref{val3}$ implies that 
\begin{equation} \label{etaorder}
\ord_{\wp_1} (\sigma \eta_1 - \tau \eta_1) \leq 2.  
\end{equation} 

Now we are in a position to compute the desired conductor.

\begin{prop} \label{localcond1}
$\cond((K_1)_{\B_1}(\widetilde{Q}_e)/(K_1)_{\B_1}) = \B_1^2$
\end{prop}

\begin{proof}
We know that $(K_1)_{\B_1}(\widetilde{Q}_e)/(K_1)_{\B_1}$ is a totally ramified extension of degree $p$.  Let $G$ represent its Galois group, $G_m$ the $m$th ramification group of $\wp_1$ over $\B_1$, and $k$ the exact order of $\B_1$ dividing $D((K_1)_{\B_1}(\widetilde{Q}_e)/(K_1)_{\B_1})$.  Then $k = \sum_{m \geq 0} (\#(G_m) -1)$.  Since $\#(G) =p$ and $\B_1$ is totally and wildly ramified, we see that $G_0 = G_1 = G$.  Therefore $k \geq 2(p-1)$ and $\B_1^{2(p-1)}$ divides $D((K_1)_{\B_1}(\widetilde{Q}_e)/(K_1)_{\B_1})$.

We know from $\eqref{etauniform}$ that $\eta_1$ generates $(K_1)_{\B_1}(\widetilde{Q}_e)$.  Let $f$ be a minimal polynomial for $\eta_1$, and $\sigma, \tau \in G$, then 
\begin{equation} \label{eq:localcond1}
\begin{split}
D((K_1)_{\B_1}(\widetilde{Q}_e)/(K_1)_{\B_1}) & =  N_{(K_1)_{\B_1}(\widetilde{Q}_e)/(K_1)_{\B_1}} (f'(\eta_1))  \\
     & = \prod_{\sigma \neq \tau } (\sigma  \eta_1 - \tau \eta_1)  \left| (\wp_1^2)^{p(p-1)} = \B_1^{2(p-1)}  \right. .  
\end{split}
\end{equation}
\noindent
Comparing equations $\eqref{eq:localcond1}$ and $\eqref{etaorder}$, we have that $\ord_{\wp_1} (\sigma \eta_1 - \tau \eta_1) = 2$, and hence $D((K_1)_{\B_1}(\widetilde{Q}_e)/(K_1)_{\B_1}) = \B_1^{2(p-1)}$. 

Now let $\psi$ be a non-trivial first degree character on $G$.  Recall that $G$ is cyclic of order $p$.  Therefore there are $p-1$ such characters, and they all have the same conductor, i.e., $\cond(\psi) = \cond((K_1)_{\B_1}(\widetilde{Q}_e)/(K_1)_{\B_1})$.  The conductor-discriminant formula and the above calculation yield:
\begin{equation}\label{K1localcond}
\B_1^{2(p-1)} = D((K_1)_{\B_1}(\widetilde{Q}_e)/(K_1)_{\B_1}) = \cond((K_1)_{\B_1}(\widetilde{Q}_e)/(K_1)_{\B_1})^{p-1}.
\end{equation}
\noindent
The result follows on taking $(p-1)$th roots of both sides of $\eqref{K1localcond}$.
\end{proof}

The following corollary is immediate, and is shown in the proof of the preceding proposition.

\begin{cor} \label{localcond2}
Let $\psi$ be a non-trivial first degree character on the Galois group of  $(K_1)_{\B_1}(\widetilde{Q}_e)/(K_1)_{\B_1}$, then $\cond(\psi) = \B_1^2$.
\end{cor}

Figure \ref{fig:localtower} begins to illustrate how we will use Proposition \ref{localcond1} and Corollary \ref{localcond2}.  The approach is to compute conductors for $n=1$, and then to use isomorphisms from field theory to lift the results to $n=e$.  It is not hard to see that we have the following isomorphism of Galois groups (for more details see \cite{Rowe}):
\begin{equation}\label{iso1}
\gal((K_e)_{\B_e}(\widetilde{Q}_e)/(K_e)_{\B_e}) \cong \gal((K_1)_{\B_1}(\widetilde{Q}_e)/(K_1)_{\B_1}),
\end{equation}
\noindent and hence
\begin{equation}\label{iso2}
\gal((K_e)_{\B_e}(\widetilde{Q}_e)/(K_1)_{\B_1}) \cong \gal((K_e)_{\B_e}(\widetilde{Q}_e)/(K_e)_{\B_e}) \times \gal((K_e)_{\B_e}/(K_1)_{\B_1}).
\end{equation}

\begin{prop} \label{condorde1}
Let $\psi$ be any non-trivial first degree character on the Galois group of $((K_e)_{\B_e}(\widetilde{Q}_e)/(K_e)_{\B_e})$ and let $\widetilde{\psi}$ be the associated character on the Galois group of $((K_1)_{\B_1}(\widetilde{Q}_e)/(K_1)_{\B_1})$ given by $\eqref{iso1}$.  Then
\begin{equation*}
\ord_{\B_e} \cond(\psi) = \ord_{\B_1} \cond(\widetilde{\psi}).
\end{equation*}
\end{prop}

\begin{proof}
Given Corollary \ref{localcond2} and Lemma \ref{Kncond}, this follows almost verbatim from the proof of Proposition 5.16 of \cite{GrantCW}.
\end{proof}

\begin{proof}[Proof of Theorem \ref{thm:cond}]
Let $\chi$ be a non-trivial first degree character of $\gal(L_e/K_e)$ such that $\chi^{p^{e-1}} \neq 1$.  By Lemma \ref{chare}, we know that $\cond(L_e/K_e) = \cond(\chi)$ and $\cond(\chi)$ is the product of its (local) Artin conductors.  Since the only ramification occurs above primes in $\SOP$, we need only compute their local conductors.  Given the assumption on $\chi$, we see that $\chi$ is non-trivial when restricted to $\gal(L_e/L_{e-1}K_e)$.  

And $\gal(L_e/L_{e-1}K_e) \cong \gal((K_e)_{\B_e}(\widetilde{Q}_e)/(K_e)_{\B_e})$, so we are in a position to apply Proposition \ref{condorde1}.   But Proposition \ref{condorde1} combined with Corollary \ref{localcond2} gives us the result.  
\end{proof}

\section{Congruence relations on units}\label{sec:units}

Let $e$ and $\SOP$ be as in section \ref{sec:filt}.  Since $e$ is fixed, we will let $F = L_e$, $E = K_e$, and $\B_i$ represent the unique totally ramified primes of $K_e$ lying over $\p_i$ of $K$ for $1 \leq i \leq g$.  Furthermore, we let $\B_1, \ldots, \B_s$ be the primes of $\SOP$.  For an ideal $\m$ of $E$, we denote the ray class field modulo $\m$ by $E(\m)$.

\begin{prop} \label{rayclass1}
$E(\B_1^2 \cdots \B_s^2)/E(\B_1 \B_2^2 \cdots \B_s^2)$ is an extension of degree $p$.
\end{prop}

\begin{proof}
By Theorem \ref{thm:cond}, we know that the extension $E(\B_1^2 \cdots \B_s^2)/E(\B_1 \B_2^2 \cdots \B_s^2)$ is non-trivial.  Indeed, $F/E$ has conductor $\cond(L/E) = \B_1^2 \cdots \B_s^2$, so we must have $E(\B_1^2 \cdots \B_s^2) \supseteq F$ and $E(\B_1 \B_2^2 \cdots \B_s^2) \not\supseteq F$.  Since $p$ is odd and of first degree, the result follows from an easy application of the snake lemma and class field theory. 
\end{proof}

\begin{cor} \label{units1}
The units of $E$ congruent to $1 \bmod \B_1 \B_2^2 \cdots \B_s^2$ are also congruent to $1 \bmod \B_1^2 \cdots \B_s^2$.
\end{cor}

\begin{rem}
Both Gupta and Grant used results similar to Corollary \ref{units1} to give congruence relations on units.  Under our hypotheses, unless $\#(\SOP)=1$, the best we can do is Theorem \ref{thm:cong}, which gives a congruence relation on an exterior product of units.
\end{rem}

\begin{lem} \label{unitscong}
Fix $i$ with $1 \leq i \leq g$ and let $v$ and $w$ be units of $E$ congruent to $1 \bmod \B_i$.  Then there exist relatively prime integers $a,b$ such that $v^a w^b \equiv 1 \bmod \B_i^2$.  
\end{lem}

\begin{proof}
By the Chinese remainder theorem, we can choose $\pi \in \OL_E$ to be a uniformizer at $\B_i$.  Let $v, w$ have the following $\B_i$-adic expansions:  $v \equiv 1 + \alpha \pi_i \bmod \B_i^2$ and $w \equiv 1 + \beta \pi_i \bmod \B_i^2$.  Since $N\B_i = p$, we may choose $\alpha, \beta \in \{ 0, \ldots, p-1 \}$.  We see that $v^a w^b \equiv 1 + (\alpha a + \beta b)\pi_i \bmod \B_i^2$.  
So by putting $a = \beta/\mbox{gcd}(\alpha,\beta)$ and $b = -\alpha/\mbox{gcd}(\alpha,\beta)$, we have $(a,b) =1$ and $v^a w^b \equiv 1  \bmod \B_i^2$.
\end{proof}

\begin{lem} \label{wedgecong}
Let $v,w$ be units of $E$ such that $v,w \equiv 1 \bmod \B_i$ and $v^a w^b \equiv 1 \bmod \B_i^2$ for some relatively prime integers $a,b$.  Then $v \wedge w = v^a w^b \wedge v^c w^d \equiv 1 \wedge v^c w^d \equiv \mathbf{1} \bmod \B_i^2$, where $c,d \in \Z$ are such that $ad - bc = 1$.  That is, $v \wedge w = v' \wedge w'$ with $v', w'$ units such that $v' \equiv 1 \bmod \B_i^2, w' \equiv 1 \bmod \B_i$.
\end{lem}

\begin{proof}
Since $(a,b)=1$, there exist integers such that $af + bg = 1$.  Let $d = f$ and $c = -g$.  By properties of exterior products, $v \wedge w = (ad-bc)(v \wedge w) = v^a w^b \wedge v^c w^d \equiv 1 \wedge v^c w^d \bmod \B_i^2$.
\end{proof}

With this result, we are now in a position to complete the prooof of the theorem.


\begin{proof}[Proof of Theorem \ref{thm:cong}]
Note that if $s=1$ then the result is trivial.  We will apply the Lemma \ref{wedgecong} repeatedly for the primes $\B_i$ with $2 \leq i \leq s$.  Starting with $\B_2$, we apply Lemma \ref{wedgecong} as many as $s-1$ times.  Then  
\begin{equation*}
u_1 \wedge \cdots \wedge u_s = u_1^{'} \wedge \cdots \wedge u_s^{'},
\end{equation*}
\noindent
where $u_1', \ldots, u_{s-1}' \equiv 1 \bmod \B_1 \B_2^2 \B_3 \cdots \B_s$ and $u_s' \equiv 1 \bmod \B_1 \cdots \B_s$.  So if we do this recursively, we get $u_1 \wedge \cdots \wedge u_s = u_1^{(s-1)} \wedge \cdots \wedge u_s^{(s-1)}$, where
\begin{equation*}
\begin{split}
u_1^{(s-1)} &  \equiv 1 \bmod \B_1\B_2^2 \cdots \B_s^2  \\
u_2^{(s-1)} &  \equiv 1 \bmod \B_1\B_2^2 \cdots \B_{s-1}^2 \B_s \\
\vdots      & \vdots \\
u_{s-1}^{(s-1)} & \equiv 1 \bmod \B_1\B_2^2 \B_3 \cdots \B_s \\
u_s^{(s-1)} & \equiv 1 \bmod \B_1\B_2 \cdots \B_s.
\end{split}
\end{equation*}
\noindent
Now by Corollary \ref{units1}, $u_1^{(s-1)} \equiv 1 \bmod \B_1^2 \cdots \B_s^2$.  Thus $u_1 \wedge \cdots \wedge u_s \equiv 1 \wedge u_2^{(s-1)} \wedge \cdots \wedge u_s^{(s-1)} \equiv \mathbf{1} \bmod \B_1^2\B_2^2\cdots\B_s^2$ as desired.
\end{proof}

\begin{cor} \label{corthm2}
Let $u_1, \ldots, u_g$ be units of $E$ congruent to $1 \bmod \B_1 \cdots \B_g $, then 
\[u_1 \wedge \cdots \wedge u_g \equiv \mathbf{1} \bmod \B_1^2 \cdots \B_g^2. \]
\end{cor}

\begin{proof}
Let $u_1, \ldots, u_g \equiv 1 \bmod \B_1 \cdots \B_g$.  By repeated application of Lemma \ref{wedgecong}, we can construct  
\begin{equation*}
u_1 \wedge \cdots \wedge u_g = \tilde{u}_1 \wedge \cdots \wedge \tilde{u}_g,
\end{equation*}
\noindent 
where $\tilde{u}_1, \ldots, \tilde{u}_s \equiv 1 \bmod \B_1 \cdots \B_s \B_{s+1}^2 \cdots \B_g^2$.  The result follows from applying Theorem \ref{thm:cong} to $\tilde{u}_1, \ldots, \tilde{u}_s$. 
\end{proof}

\bibliographystyle{amsalpha}	
\bibliography{cwrefs}		

\end{document}